\documentclass[12pt]{amsart}
\usepackage{amsfonts,amsmath,amssymb,comment,amscd,latexsym,{xy}}

\def\A{\mathbb{A}(\tl)}
\def\At{\mathbb{A}}
\def\c{\text{cc}}
\def\g{\mathcal{G}}
\def\l{\ell}
\def\gg{\mathcal{G}(\tl)}

\def\E{\mathrm{Ext}_{[\tl]}(\k C_p,\k^G)}
\def\q{\kappa}
\def\Z{\mathbb{Z}}
\def\X{\mathbb{X}(\tl)}
\def\x{\mathbb{X}}
\def\H{H^2_c(\k C_p,\k^G,\tl)}
\def\Hc{H^2_{\c}(\k C_p,\k^G,\tl)}
\def\t{\text{triv}}
\def\F{\Phi}
\def\b{\bullet}
\def\sb{{}_{\b}}
\def\we{\wedge}
\def\Tr{\text{tr}}

\newtheorem{Lem}{Lemma}[section]

\newtheorem{Prop}[Lem]{Proposition}
\newtheorem{Cor}[Lem]{Corollary}
\newtheorem{Thm}[Lem]{Theorem}

\theoremstyle{definition}
\newtheorem{Def}[Lem]{Definition}
\newtheorem{Rem}[Lem]{Remark}

\renewcommand{\k}{{\Bbbk}}
\renewcommand{\o}{\otimes}
\newcommand{\D}[1]{{\Bbbk}^{#1}}
\newcommand{\pf}{\medskip\noindent{\sc Proof:\,}}
\newcommand{\G}{\widehat G}

\newcommand{\ov}[1]{\overline {#1}}
\newcommand{\bu}[1]{{#1}^{\bullet}}

\newcommand{\bi}{\begin{itemize}}
\newcommand{\ei}{\end{itemize}}
\newcommand{\Aut}{\mathop{\mathrm{Aut}}\nolimits}
\newcommand\gen[1]{{\langle #1\rangle}}
\newcommand{\un}[1]{\underline #1}
\newcommand{\tl}{\triangleleft}
\newcommand{\tr}{\triangleright}
\numberwithin{equation}{section} 
\begin{document}

\title[Classification of Hopf algebras] {On classification of finite-dimensional semisimple Hopf algebras.}
\author{Leonid Krop}
\address{DePaul University\\Chicago, IL 60614}\email{lkrop@depaul.edu}

\begin{abstract} We develop a mechanism for classication of isomorphism types of non-trivial semisimple Hopf algebras whose group of grouplikes $G(H)$ is abelian of prime index $p$ which is the smallest prime divisor of $|G(H)|$. We describe structure of the second cohomology group of extensions of $\k C_p$ by $\k^G$ where $C_p$ is a cyclic group of order $p$ and $G$ a finite abelian group. We carry out an explicit classification for Hopf algebras of this kind of dimension $p^4$ for any odd prime $p$. The ground field is algebraically closed of characteristic $0$.
\end{abstract}
\date{3/20/15}
\maketitle



\section{Introduction}\label{intro}

We work with Hopf algebras $H$ over an algebraically closed field $\k$ of characteristic $0$. We let $G(H)$ denote the group of grouplikes of $H$. By the freeness theorem \cite{NZ} $\dim H=m\dim\k G(H)$ for an integer $m$. We say that a prime $p$ is small relative to a finite group $G$ if $p$ is the least prime divisor of $|G|$. Unless stated otherwise, we assume that $H$ is semisimple of dimension $p|G(H)|$ for a prime number $p$, $G(H)$ is abelian and $p$ is small relative to $G(H)$. For brevity, we name such Hopf algebras {\em almost abelian}. As usual, a finite-dimensional Hopf algebra is called {\em trivial} if it or its dual is a group algebra. The goal of the paper is to classify semisimple, non-trivial almost abelian Hopf algebras. 

We introduce more notation. We denote by $C_p$ a cyclic group of order $p$ and by $\k G$ and $\k^G$ the group algebra of $G$ over $\k$ and its dual, respectively. We will write $\mathrm{Ext}(\k C_p,\k^{G})$ for the set of all equivalence classes of extensions of $\k C_p$ by $\k^G$.

The problem just stated reduces to that of classifying {\em abelian extensions} of a special kind. For by a result of \cite{KM}, $\k G(H)$ is a normal subHopf algebra, a fact that combined with the theorem of Kac-Zhu \cite{{Kac},{Z}} yields that $H$ lies in $\mathrm{Ext}(\k C_p,\k^{G(H)})$. We will refer to elements of $\mathrm{Ext}(\k C_p,\k^{G(H)})$ as Hopf algebras and extensions interchangeably.

Our main concern becomes to understand the set of isomorphism types in  $\mathrm{Ext}(\k C_p,\k^G)$ where $G$ is a finite abelian group and $p$ is small relative to $G$. In general, that is for arbitrary finite groups $F,G$, there is no systematic procedure by which isomorphism classes of Hopf algebras that are extensions of $\k F$ by ${\k}^G$ can be found. One purpose of the article is to fill this gap for the case in hand. In order to formulate the statement we will require a few more notions. We write $A_p$ for the group $\mathrm{Aut}(C_p)$ of automorphisms of $C_p$. An action $\tl$ of $C_p$ on $G$ is a representation $C_p\to\mathrm{Aut}(G)$. Let $\mathcal{R}=\{\tl\}$ denote the set of all representations.  The group $\mathrm{Aut}(G)$ acts naturally on $\mathcal{R}$ by conjugation splitting $\mathcal{R}$ into the union of sets $\text{eq}(\tl)$ of representations equivalent to $\tl$. In turn, the group $A_p$ also acts on $\mathcal R$ via $\tl\mapsto\tl^{\alpha}$ where, for every $\alpha\in A_p$, $a\tl^{\alpha}x=a\tl\alpha(x),\,a\in G,x\in C_p$.
This action is passed on the sets $\text{eq}(\tl)$ via $\text{eq}(\tl)^\alpha=\text{eq}(\tl^{\alpha})$ giving rise to classes of representations $[\tl]=\cup_{\alpha}\text{eq}(\tl)^{\alpha}$. We denote the stabilizer of $\text{eq}(\tl)$ by $C(\tl)$.  
 
The splitting of $\mathcal R$ into the union of $[\tl]$ induces a splitting of

\noindent$\mathrm{Ext}(\k C_p,\k^G)$. Namely, for every $[\tl]$ we define $\E$ as the set of all equivalence classes of extensions whose $C_p$-action belongs to $[\tl]$, and then we have $\mathrm{Ext}(\k C_p,\k^G)=\displaystyle\bigcup_{[\tl]}\E$. It suffices to classify isomorphism types in each $\E$. 

To this end we bring in the second degree Hopf cohomology group \cite{{A},{M}} denoted by $\H$ after the work of M. Mastnak \cite{Ma}. We aim at constructing a subgroup $\gg$ of $\mathrm{Aut}(G)$ and its action on $\H$ compatible with isomorphism types of extensions in the sense that for any $\tau,\tau'\in\H$ $(\tau,\tl)$ and $(\tau',\tl)$ give rise to isomorphic Hopf algebras iff $\tau$ and $\tau'$ lie on the same orbit of $\gg$.

To begin with, we introduce the group $\A$ of all $C_p$-automorphisms of $(G,\tl)$. For every $\alpha\in C(\tl)$ we fix a $C_p$-isomorphism $\lambda_{\alpha}:(G,\tl)\overset{\sim}\rightarrow(G,\tl^{\alpha})$. We set $\gg$ to be the subgroup of $\mathrm{Aut}(G)$ generated by $\A$ and the set $\{\lambda_{\alpha}|\alpha\in C(\tl)\}$ if $\tl$ is nontrivial, and $\gg=\mathrm{Aut}(G)\times A_p$, otherwise.   
 
$\H$ contains a distinguished subgroup $\Hc$ of (the images of) symmetric Hopf $2$-cocycles parametrizing cocommutative extensions. Let us write ${\H/\gg}_{nc}$ for the set of $\gg$-orbits not contained in $\Hc$. Reciprocally, we let 

\noindent${}_{\text{nc}}\mathrm{Ext}_{[\tl]}(\k C_p,\k^G)/\cong$ stand for the set of isotypes of noncocommutative extensions, and we put $\text{cl}(H)$ for the isomorphism class of $H$. 

The principal result of the paper states: There is a bijection 
\begin{equation*} {\H/\gg}_{nc}\rightleftarrows  {}_{\text{nc}}\mathrm{Ext}_{[\tl]}(\k C_p,\k^G)/\cong\end{equation*}
given by  $(\tl,\tau)\g\mapsto \text{cl}(H(\tl,\tau))$. 
 
Our next concern lies with structure of $\H$. We want to find a form of $\H$ with good computational properties. We need several notions. Let $H^2(C_p,\G,\b)$ be the second degree cohomology group of extensions of $C_p$ over $\G$, where $\G$ is the dual group with the $C_p$-action $\b$ dual to $\tl$, and $H^2(G,\k^{\b})$ be the Schur multiplier of $G$. There is a mapping $N$ \cite{Mac} acting on a $\Z C_p$-module  $M$ by $N(m)=\phi_p.m$ where $\phi_p$ is the $p$th cyclotomic polynomial. We denote the kernel of $N$ in $M$ by $M_N$. The main result in the strong form states that for any {\em odd} $p$ there is a $C_p$-isomorphism
\begin{equation}\label{Hopf2cohomology}H^2_c(\k C_p,\k^G,\tl)\simeq H^2(C_p,\G,\b)\times H^2_N(G,\bu{\k})\end{equation}
We remark that this isomorphism can be seen as the Hopf cohomology version of the Baer's formula for cohomology of central extensions of a group $G$ by $C_p$ \cite [p.34]{BT}. On the other hand its main utility lies in the fact that both factors in the right-hand side of it are nicely computable in any set of generators for $G$ thanks to the classical isomorphisms $H^2(C_p,\G,\b)\simeq{\G}^{C_p}/N(\G)$  and $H^2(G,\k^{\sb})\simeq\mathrm{Alt}(G)$ \cite{Mac} and \cite{BT}, where $\mathrm{Alt}(G)$ is the group of all bilinear, alternate mapping $G\times G\to\k^{\b}$. We let $\mathbb{X}(G,\tl)$ denote the right-hand side of \eqref{Hopf2cohomology} and call it the {\em classifying group} of $\E$. For every odd $p$ $\mathbb{X}(G,\tl)$ acquires component-wise $\gg$-module structure via transport of action along the isomorphism \eqref{Hopf2cohomology}.  The new, most useful, formulation of the main theorem asserts that there is a bijection  
\begin{equation*} {\mathbb{X}(G,\tl)/\gg}_{\text{nc}}\rightleftarrows{}_{\text{nc}}{\E}/\cong\end{equation*}  
where ${\mathbb{X}(G,\tl)/\gg}_{\text{nc}}$ denotes the set orbits of $\gg$ not contained 

\noindent in ${\G}^{C_p}/N(\G)$. 

For $p=2$ less is known. We show only that the isomorphism \eqref{Hopf2cohomology} holds for elementary $2$-groups, though it is not, in general, an $\A$-isomorphism.\footnotemark\footnotetext{ See Appendix 2} 

The previous related works consist of the fundamental result of D. Stefan \cite{St} to which this article provides concrete examples, and various classification theorems. The papers \cite{{KP},{Kac},{Z},{Mas1},{Mas2},{Mas3},{Mas4}} treat a number of instances of almost abelian Hopf algebras. Namely, it is known that semisimple Hopf algebras of dimension $p$ and $p^2$ are trivial, nontrivial Hopf algebras of dimension $p^3$ are almost abelian and the number of their isomorphism types equals $p+1$ for every odd $p$, while there is a unique $8$-dimensional nontrivial Hopf algebra, and for any odd $p$ $\mathrm{Ext}(\k C_2,\k^{\mathbb{Z}_p\times\mathbb{Z}_p})$ contains a unique Hopf algebra up to isomorphism. Information on the $p^4$-dimensional semisimple Hopf algebras is limited to $p=2$ and consists of a complete classification of $16$-dimensional semisimple Hopf algebras and almost abelian Hopf algebras $H$ of dimension $2^{n+1}$ with $G(H)=\mathbb{Z}_{2^{n-1}}\times\mathbb{Z}_2$ both due to Y. Kashina \cite{{Ka},{Ka1}}. 

The paper is organized in six sections. In Section 1 we review the necessary facts of the theory of abelian extension. Sections 2 and 3 are devoted to the main results. We prove the structure theorem for the groups $H^2_c(\k C_p,\k^G,\tl)$ and the isomorphism and bijection theorems in Sections 3 and 4, respectively. Section 5 contains applications to classification of Hopf algebras of dimensions $p^2,p^3$ and an example of a non self-dual semisimple Hopf algebra of dimension $p^3$. However the bulk of this Section is devoted to finding the exact number of nontrivial almost abelian Hopf algebras of dimension $p^4$; we show that there are $5p+23$ distinct almost abelian Hopf algebras, if $p>3$, and $33$, otherwise. In the course of the proof we extend the contents of \cite{Ka1} from $p=2$ to an arbitrary prime. In the last section we revisit a theorem of Kac-Masuoka on $8$-dimensional Hopf algebras and give a generalization of a result of A. Masuoka \cite{Mas4}.
\subsection{Notation and Convention}
We adhere to the notation of \cite{Mo} on Hopf algeras and to \cite{{A},{Ma},{M}} for the theory of Hopf algebra extensions. In addition to notation in the Introduction we will use the following.

\noindent $A^{\bullet}$ the group of units of a commutative ring $A$.

\noindent $\Gamma^n$ direct product of $n$ copies of group $\Gamma$.
 
\noindent $\text{Fun}(\Gamma,A^{\bullet})$ the group of all functions from $\Gamma$ to $A^{\bullet}$ with pointwise multiplication. 
 
\noindent$Z^2(\Gamma,\bu{A},\bullet)),B^2(\Gamma,\bu{A},\bullet)$ and $H^2(\Gamma,\bu{A},\bullet)$ are the groups of $2$-cocycles, $2$-coboundaries, and the second degree cohomology group of $\Gamma$ over $\bu{A}$ with respect to an action $\bullet$ of $\Gamma$ on $A$ by ring automorphisms.

\noindent$\delta_{\Gamma}$ the differential of the standard cochain complex for cohomology of the triple $(\Gamma,\bu{A},\bullet)$ \cite[IV.5]{Mac}.

\noindent$\Z_n$ cyclic group of order $n$ additively written.

In order to simplify notation we will often use the same symbol for an element of $Z^2(\Gamma,\bu{A},\tl)$ and its image in $H^2(\Gamma,\bu{A},\tl)$. The context makes  the intended meaning clear.

Throughout the paper we treat the terms $\Gamma$-module, $\Gamma$-linear, etc as synonymous to $\Z\Gamma$-module, $\Z\Gamma$-linear, etc. We use the abbreviated term isotypes for isomophism types.
\section{Abelian extensions}\label{abelian}

In this paper we are concerned with finite-dimensional Hopf algebras over $\k$. Let $F$ and $G$ be finite groups. A Hopf algebra $H$ is an extension of $\k F$ by $\k^G$ if there is a sequence of Hopf mappings 
\begin{equation}\label{A}\k^G\overset\iota\rightarrowtail H\overset\pi\twoheadrightarrow \k F\end{equation}
with $\iota$ monomorphism, $\pi$ epimorphism, $\iota(\k^G)$ normal in $H$ and $\text{Ker}\pi=\iota(\k G)^+H$. We give a synopsis of  basic results on abelian extensions refering to \cite{M} for details. 

An abelian extension is characterized by a quadruple $D=\{\sigma,\tau,\tl,\tr\}$ called a datum for $H$ and we write $H=H(D)$. This comes about from a crossed product splitting of $H$ and $H^*$. For by \cite{NVO}, or  general theorems \cite[2.4]{S2}, \cite[3.5]{MD}\footnotemark\footnotetext{ A short independent proof is given in the Appendix 1} $H$ is a crossed product of $\k F$ over $\k^G$. Since $H^*$ is an extension of $\k G$ by $\k^F$, see \cite[4.1]{B} or \cite[3.3.1]{A}, $H^*$ is a crossed product of $\k G$ over $\k^F$. Thus there are two module algebra actions $\cdot:\k F\o\k^G\to\k^G$ and $\cdot:\k^F\o\k G\to\k^F$ and a pair of group $2$-cocycles $(\sigma,\tau)\in Z^2(F,{(\k^G)}^{\b})\times Z^2(G,{(\k^F)}^{\b})$ giving $H$ and $H^*$ an algebra structure with the multiplication
\begin{align}(f\ov{x})(f'\ov{y})&=f(x.f')\sigma(x,y))\ov{xy},\,x,y\in F,f,f'\in\k^G\label{multH}\\
(\ov{a}\phi)(\ov{b}\phi')&=\ov(ab)\tau(a,b)(\phi.b)\phi',\,a,b\in G,\phi,\phi'\in {\k}^F\label{multH^*}\end{align} 
The standard identification $\k G\cong(\k^G)^*$ via $a\mapsto\text{ev}(a):f\mapsto f(a)$ allows us to define a right action $\tl$ of $\k F$ on $\k G$ by the transpose of action $\cdot$, viz. $\gen{a\tl x,f}=\gen{\text{ev}(a),x.f}$. That is
\begin{equation}\label{dualFaction} (a\tl x)(f):=f(a\tl x)=(x.f)(a),\text{for all}\,f\in\k^G,a\in G,x\in F.\end{equation}
Likewise we obtain an action $\tr$ of $\k G$ on $\k F$. In fact both $\tl$ and $\tr$ are permutation actions on $G$ and $F$, respectively. In the dual bases $\{p_a|a\in G\}$ and $\{p_x|x\in F\}$ for $\k^G$ and $\k^F$ the two pairs of actions are related by the formulas
\begin{align}x.p_a&=p_{a\triangleleft x^{-1}}\label{.vstriangleleft}\\
p_x.a&=p_{a^{-1}\tr x}\label{.vstriangleright}.\end{align}
We fuse both actions into the definition of a product on $F\times G$ via
\begin{equation}\label{grpmult}(xa)(yb)=x(a\triangleright y)(a\triangleleft y)b\end{equation}
We use the standard notation $F\bowtie G$ for the set $F\times G$ endowed with multiplication \eqref{grpmult}.

Dualizing multiplication \eqref{multH^*} endowes $H$ with a coalgebra structure $\Delta_H,\epsilon_H$ given by \cite[4.5]{M}
\begin{align}\label{DeltaH}\Delta_H(f\ov{x})&=\sum_{a,b\in G}\tau(x,a,b)f_1p_a\ov{b\tr x}\o f_2p_b\ov{x},\\
\epsilon_H(f\ov{x})&=f(1_G).\nonumber\end{align}
We say that two structures \eqref{multH} and \eqref{DeltaH} are coherent if they turn $H$ into a bialgebra. The coherence conditions are

\begin{center}(1) $F\bowtie G$ is a group and (2) $\delta_G\sigma^{-1}=\delta_F\tau$.\end{center} 

Bialgebras so defined are always Hopf algebras, see \cite[4.7]{M} for a formula for the antipode. 

In consequence the second Hopf cohomology group of extensions \eqref{A} with fixed actions $\tl,\tr$ is defined as 
\begin{equation}\label{Hfcohomology}H^2_{\text{Hf}}(\k F,\k^G,\tl,\tr)=Z^2_{\text{Hf}}(\k F,\k^G,\tl,\tr)/B^2_{\text{Hf}}(\k f,\k^G,\tl,\tr)\end{equation}
where $Z^2_{\text{Hf}}(\k F,\k^G,\tl,\tr)=\{(\sigma,\tau)|\delta_G\sigma^{-1}=\delta_F\tau\}$ is the group of Hopf $2$-cocycles and $B^2_{\text{Hf}}(\k f,\k^G,\tl,\tr)=\{(\delta_F\zeta^{-1},\delta_G\zeta)|\zeta:F\times G\to\k^{\b}\}$ is the group of Hopf $2$-coboundaries. 

An extension \eqref{A} is called cocentral \cite{KMM} if $\k ^F$ is a central subalgebra of $H^*$. Some equivalent conditions are $\tr$ is trivial or $G$ is normal in $F\bowtie G$. Another consequence of cocentrality is that $F$ acts by Hopf automorphisms of $\k^G$ (see e.g. \cite{{Mas2},{Ka1}}). 

Our main interest lies with cocentral extensions \eqref{A} satisfying the condition
\begin{equation}\label{specialcocentral}H^2(F,\bu{(\k^G)},\tl)=\{1\}\, \text{for every action} \tl.\end{equation}
We will call them {\em special cocentral}. Below we will write $H=H(\tau,\tl)$ for a special cocentral extension with a datum $\{\tau,\tl\}$. 

In the case of special cocentral extensions the definition of cohomology groups \eqref{Hfcohomology} can be simplified. This has been done by M. Mastnak \cite{Ma} and we adopt his formulation. First we define an action of $F$ on $\mathrm{Fun}(F^n\times G^m,\k^{\b})$ extending the action $\tl$ of $F$ on $G$ via
\begin{equation}\label{Faction}y.\phi(x_1,...,x_n,a_1,...,a_m)=\phi(x_1,...,x_n,a_1\tl y,\ldots, a_m\tl y).\end{equation}
Now we let $Z^2_c(\k F,\D{G},\tl)$ and $B^2_c(\k F,\D{G},\tl)$ denote the subgroups of $Z^2(G,\bu{(\k^F)},\text{id})$ and $B^2(G,{(\k^F)}^\b,\text{id})$ of $2$-cocycles $\tau$ and $2$-coboundaries $\delta_G\eta$, respectively satisfying $\delta_F\tau=1=\delta_F\eta$.
This leads us to define  
\begin{equation*}H^2_c(\k F,\D{G},\tl)=Z^2_c(\k F,\D{G},\tl)/B^2_c(\k F,\D{G},\tl).\end{equation*}
One can see immediately that the mapping $\tau\mapsto (1,\tau)$ carries out an isomorphism between $H^2_c(\k F,\D{G},\tl)$ and $H^2_{\text{Hf}}(\k F,\k^G,\tl,\text{id})$. Explicitly both conditions $\delta_{F}\tau=\epsilon$ and $\delta_{F}\eta=\epsilon$ are expressed by:
\begin{align}\tau(xy)&=\tau(x)(x.\tau(y))\label{hopfcocycle}\\\eta(xy)&=\eta(x)(x.\eta(y))\label{hopf1cocycle}\end{align}
for all $x,y\in F$ where $F$ acts by \eqref{Faction}. The equations say that each $\tau$ and $\eta$ is a crossed homomorphism $F\to\k^{G\times G}$ and $F\to\k^G$, respectively. 

We call elements of $Z^2_c(\k F,\k^G,\tl)$ and $B^2_c(\k F,\k^G,\tl)$ {\em Hopf} $2$-cocycles and $2$-coboundaries, respectively. We will use abbreviated notation $Z^2_c(\tl),B^2_c(\tl)$, etc for $Z^2_c(\k F,\k^G,\tl)$, $B^2_c(\k F,\k^G,\tl)$, etc. when the groups $G$ and $F$ are clear from the context. We single out a subgroup $B^2_{\c}(\tl)$ of $Z^2_c(\tl)$ by the equation $B^2_{\c}(\tl)=B^2(G,(\k^F)^{\b})\cap Z^2_c(\tl)$. Clearly $B^2_c(\tl)\subset B^2_{\c}(\tl)$ so we can form the subgroup $H^2_{\c}(\tl)=B^2_{\c}(\tl)/B^2_c(\tl)$ of $H^2_c(\tl)$. We note in passing that elements of $H^2_{\c}(\tl)$ parametrize cocommutative extensions in $\mathrm{Ext}(\k F,\k^G)$.

We add a remark on $F$-invariance of subgroups just defined.
\begin{Lem}\label{Finvariance} If $F$ is abelian, then subgroups $Z^2_c(\tl),B^2_{\c}(\tl)$, and $B^2_c(\tl)$ are $F$-invariant.\end{Lem}

\pf For $Z^2_c(\tl)$ one has readily by \eqref{Faction} 
\begin{equation*}(z.\tau)(xy)=(z.\tau)(x)(zx.\tau(y)=(z.\tau)(x)(x.((z.\tau)(y))\end{equation*} 
as $x$ commutes with $z$. This shows $z.\tau\in Z^2_c(\tl)$. For the remaining two cases it suffices to note that the operator $\delta_G$ is $F$-linear on account of $G$ acting trivially on $\k^F$.\qed

\section{Structure of $H^2_c(\k C_p,\k^G,\tl)$}

From this point on $H$ is an almost abelian Hopf algebra, $G=G(H)$, $F=C_p$, and $p$ is small relative to $G$. Plainly $G$ is normal in $C_p\bowtie G$, hence the action $\tr$ is trivial. In addition, $H^2(C_p,\bu{(\k^G)},\tl)$ vanishes as $\bu{\k}$ is a divisible group by e.g. \cite[4.4]{Ma}. All in all we see that $H$ is a special cocentral extension of $C_p$ by $\k^G$.
We begin with a simple fact.
\begin{Lem} Let $\tau\in Z^2(G,\bu{(\D{C_p})})$. Then for every $x\in C_p\,\tau(x)$ is a $2$-cocycle for $G$ with coefficients in $\bu{\k}$ with the trivial action of $G$ on $\bu{\k}$.\end{Lem}
\pf The $2$-cocycle condition for the trivial action  is 
\begin{equation}\label{2cocycle} \tau(a,bc)\tau(b,c)=\tau(ab,c)\tau(a,b).\end{equation} 
Expanding both sides of the above equality in the basis $\{p_x\}$ and equating coefficients of $p_x$ proves the assertion.\qed 

Consider group $F$ acting on an abelian group $A$, written multiplicatively, by group automorphisms. Let $\mathbb{Z}F$ be the group algebra of $F$ over $\mathbb{Z}$. $\mathbb{Z}F$ acts on $A$ via
\begin{equation*} (\sum c_ix_i).a=\prod x_i.(a^{c_i}),\,c_i\in\mathbb{Z},\,x_i\in F.
\end{equation*}
For $F=C_p$ pick a generator $t$ of $C_p$ and set $\phi_i=1+t+\cdots +t^{i-1}, i=1,\ldots,p$. Choose $\tau\in Z^2(G,\bu{(\D{C_p})})$ and expand $\tau$ in 
terms of the standard basis $p_{t^i}$ for $\k^{C_p}$, $\tau=\sum\tau(t^i)p_{t^i}$ with $\tau(t^i)\in Z^2(G,\bu{\k})$. An easy induction on $i$ shows that condition \eqref{hopfcocycle} implies  
\begin{equation}\label{componentsoftau} \tau(t^i)=\phi_i.\tau(t),\;\text {for all}\;i=1,\ldots,p\end{equation}
For $i=p$ we have
\begin{equation}\label{cyclotomic} \phi_p.\tau(t)=1\end{equation}
in view of $t^p=1$ and $\tau(1)=1$.

Let $M$ be a $\mathbb{Z}C_p$-module. Following \cite{Mac} we define the mapping $N: M\to M\,\text{by}\, N(m)=\phi_p(t).m$. We denote by $M_N$ the kernel of $N$ in $M$. For $M=Z^2(G,\bu{\k}), B^2(G,\bu{\k})$ or $H^2(G,\bu{\k})$ we write $Z^2_N(G,\bu{\k})$ for $Z^2(G,\bu{\k})_N$ and similarly for the other groups. We abbreviate $Z^2_N(G,\bu{\k})$ to $Z^2_N(\tl)$ and likewise for $B^2_N(G,\bu{\k})$ and $H^2_N(G,\bu{\k})$. Thus by definition $Z^2_N(\tl)$ is the set of all $2$-cocycles satisfying  
\begin{equation}\label{Adm}\phi_p.s=1.\end{equation}
We want to compare abelian groups $Z^2_c(\tl)$ and $Z^2_N(\tl)$. This is done via the mapping 
\begin{equation*}\Theta:Z^2(G,(\k^{C_p})^{\bullet})\to Z^2(G,\k^{\bullet}),\,\Theta(\tau)=\tau(t).\end{equation*}
\begin{Lem}\label{hf=adm} The mapping $\Theta$ induces a $C_p$-isomorphism between $Z^2_c(\tl)$ and $Z^2_N(\tl)$.\end{Lem}
\pf We begin with an obvious equality $x.(\tau(y))=(x.\tau)(y)$. Taking $y=t$ we get $\Theta(x.\tau)=x.\Theta(\tau)$, that is $C_p$-linearity of $\Theta$. The relations \eqref{componentsoftau} show that $\Theta$ is monic. It remains to establish that $\Theta$ is epic.

Pick $s\in Z^2_N(\tl)$. Define $\tau: G\times G\to\bu{(\D{C_p})}$ by setting $\tau(t^i)=\phi_i(t).s,\,1\le i\le p$. The proof will be complete if we demonstrate that $\tau$ satisfies \eqref{hopfcocycle}. 

For any $i,j\le p$ we have
\begin{equation*} \tau(t^i)(t^i.\tau(t^j))=(\phi_i(t).s)(t^i\phi_j(t).s)=(\phi_i(t)+t^i\phi_j(t)).s\end{equation*}
One sees easily that $\phi_i(t)+t^i\phi_j(t)=\displaystyle\sum_{k=0}^{i+j-1}t^k$. Hence if  $i+j< p$ we have $\phi_i(t)+t^i\phi_j(t)=\phi_{i+j}(t)$ and so $\tau(t^i)(t^i.\tau(t^j))=\tau(t^{i+j})$. If $i+j=p+m$ with $m\ge 0$, then $\displaystyle\sum_{k=0}^{p+m-1}t^k=\phi_p(t)+t^p(1+\cdots+t^{m-1})$ which implies $(\displaystyle\sum_{k=0}^{p+m-1}t^k).s=\phi_p(t).s\cdot t^p\phi_m(t).s=\phi_m(t).s=\tau(t^{i+j})$ by \eqref{Adm} and as $t^p=1$.\qed

The next step is to describe structure of $H^2_{\c}(\tl)$. We need some preliminaries. First, we write $x.f$ for the left action of $C_p$ on $\k^G$ dual to $\tl$ as in \eqref{dualFaction}. Since $\G$ is the group of grouplikes of $\k^G,$ and $C_p$ acts by Hopf automorphisms $\G$ is $C_p$-stable . Further, we use
$\delta$ for the differential on the group of $1$-cochains of $G$ in $\bu{\k}$. We also note  $B^2_N(\tl)=B^2(G,\bu{\k})\cap Z^2_N(\tl)$. By \eqref{Adm} $\delta f\in B^2_N(\tl)$ iff $\phi_p(t).\delta f=1$ which, in view of $\delta$ being $C_p$-linear, is the same as $\delta(\phi_p(t).f)=1$. Since $(\delta f)(a,b)=f(a)f(b)f(ab)^{-1}$, $\text{Ker}\,\delta$ consists of characters of $G$, whence $\delta f\in B^2_N(\tl)$ iff $\phi_p(t).f$ is a character of $G$. Say $\chi=\phi_p(t).f\in \widehat{G}$. Then as $t\phi_p(t)=\phi_p(t)$, $\chi$ is a fixed point of the $C_p$-module $\widehat{G}$. Letting $\widehat{G}^{C_p}$ stand for the set of fixed points in $\G$ we have by \cite[IV.7.1]{Mac} an isomorphism 
$H^2(C_p,\widehat{G},\bullet)\simeq\widehat{G}^{C_p}/N(\widehat{G})$.  We connect $B^2_N(\tl)$ to $H^2(C_p,\widehat{G})$ via the homomorphism 
\begin{equation}\label{coboundarymap} \Phi:B^2_N(\tl)\to H^2(C_p,\widehat{G},\bullet),\,\delta f\mapsto (\phi_p.f)N(\widehat{G})\end{equation}        
\begin{Lem}\label{structureofcoboundaries} The following properties holds
\bi\item[(i)] $\Theta(B^2_{\c}(\tl))=B^2_N(\tl)$,
\item[(ii)] $\Theta(B^2_c(\tl))=\ker\Phi$,
\item[(iii)] $B^2_N(\tl)/\ker\Phi\simeq H^2(C_p,\widehat{G},\bullet)$,
\item[(iv)] $H^2_{\c}(\tl)\simeq H^2(C_p,\widehat{G},\bullet)$.
\ei \end{Lem}                                                                                                                                                
\pf First we show that $\Phi$ is well-defined. For, $\delta f=\delta g$ iff $fg^{-1}=\chi\in\widehat{G}$, hence
\begin{align*}&\Phi(\delta f)=(\phi_p.f)N(\widehat{G})=(\phi_p.g\chi)N(\widehat{G})\\
&=(\phi_p.g\cdot\phi_p.\chi)N(\widehat{G})=(\phi_p.g)N(\widehat{G})=\Phi(\delta g)\end{align*}

(i) Take some $\delta_G\eta\in B^2_{\c}(\tl)$. Evidently for every $x\in C_p$\\ (*) $(\delta_G\eta)(x)=\delta(\eta(x))$, hence $\Theta(\delta_G\eta)=\delta(\eta(t))$ is a coboundary, and $\phi_p.\delta(\eta(t))=1$ by \eqref{cyclotomic}, whence $\Theta(\delta_G\eta)\in B^2_N(\tl)$.
Conversely, pick $\delta f\in B^2_N(\tl)$ and define $\omega=\sum_{i=1}^p(\phi_i.\delta f)p_{t^i}$. The argument of Lemma \ref{hf=adm} shows $\omega$ lies in $Z^2_c(\tl)$. Set $\eta=\sum_{i=1}^p(\phi_i.f)p_{t^i}$. Using (*) again we derive 
\begin{equation*}\delta_G\eta=\sum_{i=1}^p(\phi_i.\delta f)p_{t^i}=\omega, \end{equation*}
hence $\delta_G\eta\in B^2_{\c}(\tl)$. Clearly $\Theta(\delta_G\eta)=\delta f$.

(ii) The argument of Lemma \ref{hf=adm} is applicable to $1$-cocycles satisfying \eqref{hopf1cocycle}. It shows that $\eta$ satisfies \eqref{hopf1cocycle} iff
\begin{equation}\label{hopf1cocycle'}\eta(t^i)=\phi_i.\eta(t)\end{equation}                                                                                  For $i=p$ we get $\phi_p.\eta(t)=\epsilon$, hence the calculation
\begin{align*}\Phi(\Theta(\delta_G\eta))=\Phi(\delta(\eta(t)))=(\phi_p.\eta(t)) N(\widehat{G})=N(\widehat{G}).\end{align*}
gives one direction. Conversely, $\Phi(\delta f)\in N(\widehat{G})$ means $\phi_p.f=\phi_p.\chi$ which implies $\phi_p.f\chi^{-1}=\epsilon$. Set $g=f\chi^{-1}$ and define   $1$-cocycle $\eta_g=\sum_{i=1}^p(\phi_i.g)p_{t^i}$. Since $\phi_p.g=\epsilon$,  $\eta_g$ satisfies \eqref{hopf1cocycle}, whence  $\delta_G\eta_g\in B^2_c(\tl)$. As $(\delta_G\eta_g)(t)=\delta g=\delta f$ by construction, $\Theta(\delta_G\eta_g)=\delta f$.                                                                                   
            
(iii) We must show that $\Phi$ is onto. For every character $\chi$ in $\widehat{G}^{C_p}$ we want to construct an $f:G\to\bu{\k}$ satisfying $\phi_p.f=\chi$. To this end we consider splitting of $G$ into the orbits under the action of $C_p$. Since every orbit is either regular, or a fixed point we have
\begin{equation*} G=\cup_{i=1}^r\{g_i,g_i\tl t,\ldots,g_i\tl t^{p-1}\}\cup G^{C_p}\end{equation*} 
For every $s\in G^{C_p}$ we pick a $\rho_s\in\k$ satisfying $\rho^p_s=\chi(s)$. We define $f$ by the rule
\begin{align*}f(g_i)=\chi(g_i),\,f(g_i\tl t^j)&=1\;\text{for all}\;j\ne 1\,\text{and all}\;i=1,\ldots,r,\,\text{and}\\
f(s)&=\rho_s\;\text{for every}\;s\in G^{C_p}\end{align*}                                                                                                     By definition $(\phi_p.f)(g)=\prod_{j=0}^{p-1}f(g\tl t^j)$. Therefore $(\phi_p.f)(s)=\rho_s^p=\chi(s)$ for every $s\in G^{C_p}$. If $g=g_i\tl t^j$ for some $i,j$, then a calculation $(\phi_p.f)(g)=f(g_i)=\chi(g_i)=\chi(g_i\tl t^j)=\chi(g)$, which uses the fact that $\chi$ is a fixed point under the action by $C_p$, completes the proof.

(iv) follows immediately from $H^2_{\c}(\tl)=B^2_{\c}/B^2_c(\tl)$ and parts (i)-(iii).\qed
\begin{Cor}\label{basicisic} Isomorphism $\Theta$ induces a $C_p$-isomorphism 

\noindent$\Theta_*:H^2_c(\tl)\simeq Z^2_N(\tl)/{\ker\Phi}.$\end{Cor}\qed
 
We proceed to the main result of the section. 
\begin {Prop}\label{hopf2cohomology} Suppose $G$ is a finite abelian group. If $|G|$ is odd, or $G$  is a $2$-group and either $C_2$-action is trivial, or $G$ is an elementary $2$-group, there exists a $C_p$-isomorphism
\begin{equation}\label{structureofcohomology}H^2_c(\tl)\simeq H^2(C_p,\G,\bullet)\times H^2_N(G,\bu{\k}).\end{equation}
\end{Prop}
\pf (1) First we take up the odd case. By the preceeding Corollary we need to decompose $Z^2_N(\tl)/{\ker\Phi}$. We note that for any $p$ and $G$ there is a group splitting $Z^2(G,\bu{\k})=B^2(G,\bu{\k})\times H^2(G,\bu{\k})$ due to the fact that the group of $1$-cocycles ${\bu{\k}}^{G}$ is injective, and hence so is $B^2(G,\bu{\k})$. We aim at finding a $C_p$-invariant complement to $B^2(G,\bu{\k})$. To this end we recall a well-known isomorphism $\un{a}:H^2(G,\bu{\k})\tilde{\to}\text{Alt}(G)$, see e.g. \cite[\S 2.3]{Y}. There $\text{Alt}\,(G)$ is the group of all bimultiplicative alternating functions 
\begin{equation*}\beta:G\times G\to\bu{\k},\,\beta(ab,c)=\beta(a,c)\beta(b,c),\text{and}\,\beta(a,a)=1\,\text{for all}\,a\in G.\end{equation*}
For the future applications we outline the construction of $\un{a}$. Namely, $\un{a}$ is the antisymmetrization mapping sending $z\in Z^2(G,\bu{\k})$ to $\un{a}(z)$ defined by $\un{a}(z)(a,b)=z(a,b)z^{-1}(b,a)$. One can check that $\un{a}(z)$ is bimultiplicative (cf. \cite[(10)]{Y}) and it is immediate that $\un{a}$ is $C_p$-linear. Another verification gives $\text{im}\,\un{a}=\text{Alt}(G)$ and, moreover, 

\noindent$\ker\,\un{a}=B^2(G,\bu{\k})$, see \cite[Thm.2.2]{Y}.  Thus we obtain a $C_p$-isomorphism $H^2(G,\bu{\k})\simeq\text{Alt}(G)$.

Since elements of $\text{Alt}(G)$ are bimultiplicative mappings $\text{Alt}(G)\subset Z^2(G,\k^{\b})$. For every $\beta\in\text{Alt}(G)$ a simple calculation gives $\un{a}(\beta)=\beta^2$. Thus $\un{a}(\beta)\ne 1$ as the order of $\beta$ divides the exponent of $G$. It follows $B^2(G,\bu{\k})\cap\text{Alt}(G)=\{1\}$ which gives a splitting of abelian groups
\begin{equation*}Z^2(G,\bu{\k})=B^2(G,\bu{\k})\times \text{Alt}(G)\end{equation*}
But now both subgroups $B^2(G,\bu{\k})$ and $\text{Alt}(G)$ are $C_p$-invariant hence there holds $Z^2_N(G,\bu{\k})=B^2_N(G,\bu{\k})\times \text{Alt}_N(G)$ which, in view of\\ $\text{Alt}(G)=H^2(G,\bu{\k})$, is the same as
\begin{equation}\label{Z^2adm} Z^2_N(\tl)=B^2_N(\tl)\times H^2_N(G,\bu{\k}).\end{equation} 
Now part (iii) of Lemma \ref{structureofcoboundaries} completes the proof of (1).

(2) Here we prove the second claim of the Proposition. We decompose $G$ into a product of cyclic groups $\gen{x_i},1\le i\le m$. For every $\alpha\in\text{Alt}(G)$ we define $s_{\alpha}\in Z^2(G,\k^{\bullet})$ via
$$s_{\alpha}(x_i,x_j)=\begin{cases}\alpha(x_i,x_j),&\text{if}\;i\le j\\\phantom{\alpha}1,&\text{else}.\end{cases}$$
Since $s_{\alpha}\cdot s_{\beta}=s_{\alpha\beta}$ the set $S=\{s_{\alpha}|\alpha\in\text{Alt}(G)\}$ is a subgroup of $Z^2(G,\k^{\bullet})$. One can see easily that $s_{\alpha}=s_{\beta}\Leftrightarrow\alpha=\beta$ and $\un{a}(s_{\alpha})=\alpha$, hence $S$ is isomorphic to $\text{Alt}(G)$ under $\un{a}$. For every $z\in Z^2_N(\t),\un{a}(z)\in\text{Alt}_N(G)$, and therefore $\un{a}(z)=\un{a}(s)$ for some $s\in S_N$. We have $zs^{-1}\in B^2(G,\k^{\bullet})$, but as $zs^{-1}$ has order $2$, $zs^{-1}\in B^2_N(\t)$. Thus $Z^2_N(\t)=B^2_N(\t)\times S_N$ which proves \eqref{structureofcohomology}.

(3) We prove the last claim of the Proposition. Below $G$ is an elementary $2$-group, and action of $C_2$ is nontrivial. First we establish an intermediate result, namely
\begin{Lem}\label{nonsplitcase} If action $\tl$ is nontrivial, then $Z^2_N(\tl)$ is a nonsplit extension of $\text{Alt}_N(G)$ by $B^2_N(\tl)$.\end{Lem}
\pf This will be carried out in steps.

(i) We aim at finding a basis for $\text{Alt}_N(G)$. We begin by noting that as $\text{Alt}(G)$ has exponent $2$, $\text{Alt}_N(G)$ is the set of all fixed points in $\text{Alt}(G)$. Put $R=\mathbb{Z}_2C_2$. One can see easily that $R$-module $G$ decomposes as
\begin{equation}\label{decomposition}G=R_1\times\cdots\times R_m\times G_0\end{equation}
where  $R_i\simeq R$ as a right $C_2$-module, and $G_0=G^{C_2}$. Denote by $t$ the generator of $C_2$. For each $i$ let $\{x_{2i-1},x_{2i}\}$ be a basis of $R_i$ such that $x_{2i-1}\tl t=x_{2i}$. We also fix a basis $\{x_{2m+1},\ldots,x_n\}$ of $G_0$.

We associate to every subset $\{i,j\}$ the bilinear form $\alpha_{ij}$ by setting
\begin{equation*}\alpha_{ij}(x_i,x_j)=\alpha_{ij}(x_j,x_i)=-1,\,\text{and}\,\alpha_{ij}(x_k,x_l)=1\,\text{for any}\,\{k,l\}\!\ne\{i,j\}.\end{equation*}
The set $\{\alpha_{ij}\}$ forms a basis of $\text{Alt}(G)$. One can check easily that $t$ acts on basic elements as follows
\begin {equation}\label{actioninAlt} t.\alpha_{ij}=\alpha_{kl}\;\text {if and only if}\,\{x_i,x_j\}\tl t:=\{x_i\tl t,x_j\tl t\}=\{x_k,x_l\}.\end{equation}
Recall the element $\phi_2=1+t\in\mathbb{Z}C_2$. We define forms $\beta_{ij}$ via
\begin{equation}\label{basisAlt}\beta_{ij}=\phi_2.\alpha_{ij}\;\text{if}\;t.\alpha_{ij}\ne\alpha_{ij},\;\text{and}\;\beta_{ij}=\alpha_{ij},\;\text{otherwise}.\end{equation} 
The label $ij$ on $\beta_{ij}$ is not unique as $\beta_{ij}=\beta_{kl}$ whenever $\{x_i,x_j\}\tl t=\{x_k,x_l\}$. Of the two sets $\{i,j\}$ and $\{k,l\}$ labeling $\beta_{ij}$ we agree to use the one with the smallest element, and call such minimal. We claim:
\begin{equation}\label{basisALTN} \text{The elements}\; \{\beta_{ij}\}\; \text{form a basis of}\; \text{Alt}_N(G).\end{equation}
\pf  First we note that for every group $M$ of exponent $2$ $M_N=M^{C_2}$. Suppose $\beta\in\text{Alt}(G)^{C_2}$. Say $\beta=\prod\alpha_{ij}^{e_{ij}},\;e_{ij}=0,1$. From $t.\beta=\prod(t.\alpha_{ij})^{e_{ij}}=\beta$ we see that if $\alpha_{ij}$ occurs in $\beta$, i.e. $e_{ij}=1$, then so does $t.\alpha_{ij}$, hence $\beta$ is a product of $\beta_{ij}$.\qed

(ii) We want to show $\un{a}(Z^2_N(\tl))=\text{Alt}_N(G)$. The restriction $\un{a}^*$ of $\un{a}$ to $Z^2_N(\tl)$ induces a $C_2$-homomorphism $Z^2_N(\tl)\overset{\un{a}^*}\rightarrow\text{Alt}_N(G)$ whose kernel equals $B^2(G,\k^{\bullet})\cap Z^2_N(\tl)=:B^2_N(\tl)$. 

We begin by showing $\phi_2.\text{Alt}(G)\subset \text{im}\,\un{a}^*$. For, if $\beta=\phi_2.\alpha$, pick an $s\in Z^2(G,\k^{\bullet})$ with $\un{a}(s)=\alpha$. Then $(t-1).s\in Z^2_N(\tl)$, and $\un{a}((t-1).s)=(t-1).\un{a}(s)=(t-1).\alpha=\phi_2.\alpha$, as $\alpha^2=1$, which gives the inclusion.

By step (i) and definition \eqref{basisAlt} it remains to show that all fixed points $\alpha_{ij}$ lie in $\text{im}\,\un{a}^*$. By formula \eqref{actioninAlt} $\alpha_{ij}$ is a fixed point if and only if either 
\begin{equation*}(a)\quad \{i,j\}\subset\{2m+1,\ldots,n\}\;\text{or}\; (b)\quad \{i,j\}=\{2k-1,2k\}\end{equation*} 
for some $k,\,1\le k\le m$. Below we find it convenient to write $s_{i,j}$ for $s_{\alpha_{ij}}$.

Consider case (a). We claim $s_{i,j}$ is a fixed point. For, $t.s_{i,j}$ is bimultiplicative, hence is determined by its values at $(x_k,x_l)$. It is immediate that $t.s_{i,j}(x_k,x_l)=s_{i,j}(x_k,x_l)$ for all $(x_k,x_l)$, whence the assertion. Since $s_{i,j}^2=1$ for all $i,j$, $\phi_2.s_{i,j}=1$, hence $s_{i,j}\in Z^2_N(\tl)$. As $\un{a}(s_{i,j})=\alpha_{ij}$, this case is done.

We take up (b). Say $z=s_{2i-1,2i}$ for some $i, 1\le i\le m$. An easy verification gives $\phi_2.z=\alpha_{2i-12i}\ne 1$. Thus $z\notin Z^2_N(\tl)$. To prove (ii) we need to find a coboundary $\delta g_i$ such that $z\delta g_i\in Z^2_N(\tl)$. Since $\un{a}(\alpha_{2i-12i})=1,\alpha_{2i-12i}=\delta f_i$ for some $f_i:G\to\bu{\k}$. Put $G_i$ for the subgroup of $G$ generated by all $x_j,j\ne 2i-1,2i$. We assert that one choice is the function $f_i$ defined by
\begin{equation}\label{specialcocycle} f_i(x_{2i-1}^{j_1}x_{2i}^{j_2}x')=(-1)^{j_1+j_2+j_1j_2}\,\text{for all}\;x'\in G_i\end{equation}
For, on the one hand it is immediate that for any $x',x''\in G_i$ \begin{equation*}\alpha_{2i-12i}(x_{2i-1}^{j_1}x_{2i}^{j_2}x',x_{2i-1}^{k_1}x_{2i}^{k_2}x'')=(-1)^{j_1k_2+j_2k_1}\end{equation*}
On the other hand the definitions of $f_i$ and differential $\delta$ give
\begin{eqnarray*}\lefteqn{\delta  f_i(x_{2i-1}^{j_1}x_{2i}^{j_2}x',x_{2i-1}^{k_1}x_{2i}^{k_2}x'')}\\&&=(-1)^{j_1+j_2+j_1j_2}(-1)^{k_1+k_2+k_1k_2}(-1)^{j_1+k_1+j_2+k_2+(j_1+k_1)(j_2+k_2)}\\&&=(-1)^{j_1k_2+j_2k_1}\nonumber\end{eqnarray*}
Define the function $g_i:G\to\bu{\k}$ by $g_i(x_{2i-1}^{j_1}x_{2i}^{j_2}x')=\iota^{j_1+j_2+j_1j_2}$ where $\iota^2=-1$. One can check easily the equalities $f_i^2=1$ and $t.g_i=g_i,\,g_i^2=f_i$. Hence we have $f_i(\phi_2.g_i)=f_ig_i^2=f_i^2=1$, and then a calculation 
\begin{equation*} \phi_2.(z\delta g_i)=(\phi_2.z)(\phi_2.\delta g_i)=\delta f_i\cdot\delta(\phi_2.g_i)=\delta(f_i(\phi_2.g_i))=1\end{equation*}
completes the proof of (ii).

(iii) Suppose $Z^2_N(\tl)=B^2_N(\tl)\times C$ where $C$ is a $C_2$-invariant subgroup. Then $C$ is mapped isomorphically on $\text{Alt}_N(G)$ under $\un{a}$ and so there is a unique $z\in C$ such that $\un{a}(z)=\alpha_{12}$. Since $\un{a}(s_{1,2})=\alpha_{12},\,z=s_{1,2}\delta g$ for some $g:G\to\bu{\k}$. Further, as $\alpha_{12}$ is a fixed point $\un{a}(t.z)=\alpha_{12}$ as well, hence $t.z=z$. In addition, since $\text{Alt}(G)$ is an elementary $2$-group, 
$1=z^2=(s_{1,2}\delta g)^2=(\delta g)^2=\delta(g^2)$. It follows that $g^2$ is a character of $G$. Moreover, $t.z=z$ is equivalent to $t.s_{1,2}(t.\delta g)=s_{1,2}\delta g$ which in turn gives $s_{1,2}(t.s_{1,2})(t.\delta g)=\delta g$. As $\phi_2.s_{1,2}=\alpha_{12}=\delta f_1$ we have $\delta f_1(t.\delta g)=\delta g$ which implies $\delta f_1=\delta g(t.\delta g)$ on the account of $(\delta g)^2=\delta(g^2)=1$ as $g^2$ is a character. Equivalently we have the equality 
\begin{equation}\label{splitting} f_1=g\cdot(t.g)\cdot\chi\;\text{for some}\;\chi\in\G.\end{equation}
Noting that $f_1$ is defined up to a character of $G$ we can assume that $f_1(x_1)=1=f_1(x_2)\;\text{and}\;f_1(x_1x_2)=-1$. For, $f_1$ is defined as any function satisfying $\delta f_1=\alpha_{12}$. As $\delta(f_1\chi)=\delta f_1$ for any $\chi\in\G$, $f_1$ can be modified by any $\chi$. By \eqref{specialcocycle} $f_1(x_j)=-1=f_1(x_1x_2),j=1,2$ so we can take $\chi$ such that $\chi(x_1)=\chi(x_2)=-1$. The equality \eqref{splitting} implies that for some $\chi\in\G$ there holds  
\begin{align}1&=f_1(x_j)=g(x_1)g(x_2)\chi(x_j),\,j=1,2\tag{*},\, \text{and}\\ 
-1&=f_1(x_1x_2)=g(x_1x_2)^2\chi(x_1x_2)\tag{**}\end{align} 
as $t$ swaps $x_1$ and $x_2$. Since $g^2$ is a character, $g^2(a)=\pm 1$ for every $a\in G$. It follows that $g(x_1)=\iota^m$ and $g(x_2)=\iota^k$ for some $0\le m,k\le 3$. Then equation (*) gives $1=\iota^{m+k}\chi(x_j)$. This equality shows that $\chi(x_1)=\chi_(x_2)$ and $m+k$ is even, because $\chi(a)=\pm 1$ for all $a$. Now (**), and the fact that $g^2$ is a character, gives $-1=g^2(x_1)g^2(x_2)\chi(x_1)\chi(x_2)=\iota^{2(m+k)}\iota^{-2(m+k)}=1$, a contradiction. This completes the proof of the Lemma.\qed

Finally we prove (3). Let $G$ be a group with a decomposition \eqref{decomposition}. Set $C$ to be the subgroup of $Z^2_N(\tl)$ generated by the set $B=B'\cup B''\cup B'''$ where
\begin{align*}B'&=\{\phi_2.s_{i,j}|\alpha_{ij}\,\text{is not a fixed point, and}\,\{i,j\}\,\text{is minimal}\} \\B''&=\{s_{i,j}|i<j\,\text{and}\,\{i,j\}\subset\{2m+1,\ldots,n\}\}\\B'''&=\{s_{2i-1,2i}\delta g_i|i=1,\ldots,m\}.\end{align*} 
There $g_i$ is chosen as in the case (ii) of Lemma \ref{nonsplitcase}. Passing on to $Z^2_N(\tl)/\ker\Phi$ we denote by $\ov{B^2_N(\tl)}$ and $\ov{C}$ the images of these subgroups in $Z^2_N(\tl)/\ker\Phi$. Pick a $v\in B$. If $v\in B'\cup B''$ then $v^2=1$ because the corresponding $s_{i,j}$ has order $2$. For $v=s_{2i-1,2i}\delta g_i$, $v^2=\delta g_i^2=\delta f_i$. We know $t.f_i=f_i$ and $f_i^2=1$ and therefore $\phi_2.f_i=1$, whence $\delta f_i\in\ker\Phi$ by definition \eqref{coboundarymap}. It follows that ${\ov{v}}^2=1$ for all $\ov{v}\in\ov{B}$. Furthermore, by Lemma \ref{nonsplitcase} the mapping $\un{a}$ sends $\ov{B}$ to the basis \eqref{basisALTN} of $\text{Alt}_N(G)$. Therefore $\ov{C}$ is isomorphic to $\text{Alt}_N(G)$ at least as an abelian group and forms a complement to $\ov{B^2_N(\tl)}$ in $Z^2_N(\tl)/\ker\Phi$. Since $\text{Alt}_N(G)$ consists of fixed points the proof will be completed if we show the same for $\ov{C}$. The fact that $B'\cup B''$ consists of fixed points follows from $t\phi_2=\phi_2$ and part (a) of Lemma \ref{nonsplitcase}(ii). For an $s_{2i-1,2i}\delta g_i$, the equality $\phi_2.s_{2i-1,2i}=\delta f_i$ gives $t.s_{2i-1,2i}=s_{2i-1,2i}\delta f_i$. Since $\delta f_i\in\ker\Phi$ and $t.\delta g_i=\delta g_i$ we see that $\overline{s_{2i-1,2i}\delta g_i}$ is a fixed point in $Z^2_N(\tl)/\ker\Phi$ which completes the proof.\qed 
                                                                                                                                           
\section{The Isomorphism Theorems}\label{Main}
We begin with a general observation. Let $H$ be be an extension of type (A). The mapping $\pi$ induces a $\k F$-comodule structure $\rho_{\pi}$ on $H$ via
\begin{equation}\label{comodulealg}\rho_{\pi}:H\to H\o\k F,\,\rho_{\pi}(h)=h_1\o\pi(h_2).\end{equation}
$H$ becomes an $F$-graded algebra with the graded components $H_f=\{h\in H|\rho_{\pi}(h)=h\o f\}$. Let $\chi:\k F\to H$ be a section of $\k F$ in $H$. By definition $\chi$ is a convolution invertible $\k F$-comodule mapping, that is
\begin{equation}\label{Section}\rho_{\pi}(\chi(f))=\chi(f)\o f,\,\text{for every}\,f\in F\end{equation}
Set $\ov{f}=\chi(f)$. The next lemma is similar to \cite[3.4]{M} or \cite[7.3.4]{Mo}.
\begin{Lem}\label{gradedcomponents} For every $f\in F$ there holds $H_f=\D{G}\ov{f}$\end{Lem}
\pf By definition of components $H_1=H^{\text{co}\pi}$ which equals to $\D{G}$ by the definition of extension. By \eqref{Section} $\rho_{\pi}(\ov{f})=\ov{f}\o f$, hence $\D{G}\ov{f}\subset H_f$. Since the containment holds for all $f$, the equalities 
\begin{equation*}H=\oplus_{f\in F}H_f=\oplus_{f\in F}\D{G}\ov{f}\end{equation*}
force $H_f=\D{G}\ov{f}$ for all $f\in F$.\qed
\begin{Def}Given two $F$-graded algebras $H=\oplus H_f$ and $H'=\oplus H'_f$ and an automorphism $\alpha:F\to F$ we say that a linear mapping $\psi:H\to H'$ is an $\alpha$-graded morphism if $\psi(H_f)=H'_{\alpha(f)}$ for all $f\in F$.\end{Def}
\begin{Lem}\label{gradedpsi} Suppose $H$ and $H'$ are two extensions of $\k F$ by $\D{G}$ and $\psi:H\to H'$ a Hopf isomorphism sending $\D{G}$ to $\D{G}$. Then $\psi$ is an $\alpha$- graded mapping for some $\alpha$.\end{Lem}
\pf Suppose $H$ and $H'$ are given by sequences
\begin{equation*}\k^G\overset\iota\rightarrowtail H\overset\pi\twoheadrightarrow \k F,\,\text{and}\quad\k^G\overset{\iota'}\rightarrowtail H'\overset{\pi'}\twoheadrightarrow \k F\end{equation*}
By definition of extension $\text{Ker}\,{\pi}=H(\D{G})^+$ and likewise $\text{Ker}\,{\pi'}=H'(\D{G})^+$. By assumption $\psi(\D{G})=\D{G}$, hence $\psi$ induces a Hopf isomorphism $\alpha:H/H(\D{G})^+\to H'/H'(\D{G})^+$. Replacing $H/H(\D{G})^+$ and $ H'/H'(\D{G})^+$ by $\k F$ we can treat $\alpha$ as a Hopf isomorphism $\alpha:\k F\to\k F$. $\alpha$ is in fact an automorphism of $F$. We arrive at a commutative diagram
\begin{equation*}\begin{CD}
{\k}^G @>\iota>> H @>\pi>>\k F\\
@V{\psi}VV @V{\psi}VV  @V{\alpha}VV \\      
{\k}^G @>\iota'>> H' @>\pi'>>\k F
\end{CD}\end{equation*}
Since $\psi$ is a coalgebra mapping for every $f\in F$ we have 
\begin{align*}\Delta_{H'}(\psi(\ov{f}))&=(\psi\o\psi)\Delta_H(\ov{f})=\psi((\ov{f})_1)\o\psi((\ov{f})_2),\,\text{hence}\\
\rho_{\pi'}(\psi(\ov{f}))&=\psi((\ov{f})_1)\o \pi'\psi((\ov{f})_2)=\psi((\ov{f})_1)\o\alpha\pi((\ov{f})_2)\end{align*}
On the other hand, applying $\psi\o\alpha$ to the equality 

\noindent$\rho_{\pi}(\ov{f})=(\ov{f})_1\o\pi((\ov{f})_2)=\ov{f}\o f$ gives
\begin{equation*}\psi((\ov{f})_1)\o\alpha\pi((\ov{f})_2)=\psi(\ov{f})\o\alpha(f)\end{equation*}
whence we deduce $\rho_{\pi'}(\psi(\ov{f}))=\psi(\ov{f})\o\alpha(f)$. Thus $\psi(\ov{f})\in H'_{\alpha(f)}$ which shows the inclusion 
\begin{equation*}\psi(H_f)=\psi(\D{G}\ov{f})=\D{G}\psi(\ov{f})\subseteq H'_{\alpha(f)}=\D{G}\ov{\alpha(f)}\end{equation*}
Since both sides of the above inclusion have equal dimensions, the proof is complete.\qed

In what follows $H$ is an almost abelian Hopf algebra, $G=G(H)$, $F=C_p$, and $p$ is small relative to $G$. Let $\tl$ and $\tl'$ be two actions of $C_p$ on $G$. We denote $(G,\tl)$ and $(G,\tl')$ the corresponding $C_p$-modules and we use the notation `$\bullet$' and `$\circ$' for the actions of $C_p$ on $\D{G}$ corresponding by \eqref{dualFaction} to $\tl$ and $\tl'$, respectively. We let $I(\tl,\tl')$ denote the set of all automorphisms of $G$ intertwining actions $\tl$ and $\tl'$, that is automorphisms $\lambda:G\to G$ satisfying
\begin{equation}\label{modulemap}(a\tl x)\lambda=a\lambda\tl' x,\,a\in G,x\in C_p\end{equation}
We make every $\lambda$ act on  functions $\tau:C_p\times G^2\to\k^{\bullet}$ by 
\begin{equation*}(\tau.\lambda)(x,a,b)=\tau(x,a\lambda^{-1},b\lambda^{-1}).\end{equation*} 

\begin{Lem}\label{actiononcocycles} {\rm(i)} The group $Z^2(G,\bu{(\D{C_p})})$ is invariant under the action induced by any automorphism of $G$,

{\rm(ii)} A $C_p$-isomorphism $\lambda:(G,\tl)\to (G,\tl')$ induces $C_p$-isomorphisms between the groups $Z^2_c(\tl),B^2_c(\tl),H^2_c(\tl)$ and $Z^2_c(\tl'),B^2_c(\tl'),H^2_c(\tl')$, respectively.
\end{Lem}
\pf (i) is immediate.

(ii) We must check condition \eqref{hopfcocycle} for $\tau.\lambda$ and $\mathbb{Z}C_p$-linearity of the induced map. First we note $\lambda^{-1}$ is a $C_p$- isomorphism between $(G,\tl')$ and $(G,\tl)$, as one can check readily. Next we verify \eqref{hopfcocycle} and $C_p$-linearity in a single calculation
\begin{align*} &(\tau.\lambda)(xy)(a,b)=\tau(xy,a\lambda^{-1},b\lambda^{-1})\\
&=\tau(x,a\lambda^{-1},b\lambda^{-1})(x\bullet\tau(y,a\lambda^{-1},b\lambda^{-1}))\\
&=\tau(x,a\lambda^{-1},b\lambda^{-1})\tau(y,a\lambda^{-1}\tl x,b\lambda^{-1}\tl x)\\&=\tau(x,a\lambda^{-1},b\lambda^{-1})\tau(y,(a\tl' x)\lambda^{-1},(b\tl' x)\lambda^{-1})\\
&=(\tau.\lambda)(x)(x\circ(\tau.\lambda)(y))(a,b).
\end{align*}                                                                                                                                                 In the case of $B^2_c(\tl)$, first one checks the equality
\begin{equation*}(\delta_{G}\eta).\lambda=\delta_{G}(\eta.\lambda)\,\text{for any}\,\eta:C_p\times G\to\bu{\k}.\end{equation*}                                                        
It remains to verify the condition \eqref{hopf1cocycle} for $\eta.\lambda$. That is done similarly to the calculation in (ii).\qed

Let $(G,\tl)$ be a $C_p$-module. Recall that $\mathbb{A}(\tl)$ denotes the group of $C_p$-automorphisms of $(G,\tl)$. By the above Lemma $Z^2_c(\tl)$ is an $\mathbb{A}(\tl)$-module. Symmetrically, the group $A_p=\mathrm{Aut}(C_p)$ of automorphisms of $C_p$ acts on $\mathrm{Map}(C_p\times G^2,\bu{\k})$ via
\begin{equation*}\tau.\alpha(x,a,b)=\tau(\alpha(x),a,b)\end{equation*}
We want to know the effect of this action on $Z^2_c(\tl)$. Let $(G,\tl)$ be a $C_p$-module. For $\alpha\in A_p$ we define a $C_p$-module $(G,\tl^{\alpha})$ via
\begin{equation*}a\tl^{\alpha}x=a\tl\alpha(x),\;a\in G,\,x\in C_p\end{equation*}
Similarly, an action `$\bullet$' of $C_p$ on $\k^G$ can be twisted by $\alpha$ into `$\bullet^{\alpha}$' via 
\begin{equation*}x\bullet^{\alpha}r=\alpha(x)\bullet r,\;r\in\k^G\end{equation*}  
One can see easily that if $\bullet$ and $\tl$ correspond to each other by \eqref{dualFaction}, then so do $\bullet^{\alpha}$ and $\tl^{\alpha}$. 
\begin{Lem}\label{secondaction} {\rm(i)} If $\lambda\in I(\tl,\tl')$, then $\lambda\in I(\tl^{\alpha},\tl'^{\alpha})$ for every $\alpha\in A_p$,
 
\rm{(ii)} The mapping $\tau\mapsto\tau.\alpha$ induces an $\mathbb{A}(\tl)$-isomorphism between $Z^2_c(\tl),B^2_c(\tl),H^2_c(\tl)$ and $Z^2_c(\tl^{\alpha}),B^2_c(\tl^{\alpha}),H^2_c(\tl^{\alpha})$, respectively for every $\alpha\in A_p$ .\end{Lem}
\pf {\rm(i)} For every $a\in G,x\in C_p$ we have
\begin{equation*}(a\tl^{\alpha})\lambda=(a\tl\alpha(x))\lambda=a\lambda\tl'\alpha(x)=a\lambda\tl'^{\alpha}x\end{equation*}

{\rm(ii)} First we note that $\mathbb{A}(\tl)$ can be identified with $\mathbb{A}(\tl^{\alpha})$ for any $\alpha$ by the folllowing calculation
\begin{equation*}(g\tl^{\alpha}x)\phi=(g\tl\alpha(x))\phi=(g\phi)\tl\alpha(x)=g\phi\tl^{\alpha}x\,\text{for every}\,\phi\in\mathbb{A}(\tl).\end{equation*}
Thus we will treat every $Z^2_c(\tl^{\alpha})$ as an $\mathbb{A}(\tl)$- module. Our next step is to show that for every $\tau\in Z^2_c(\tl)$, $\tau.\alpha$ lies in $Z^2_c(\tl^{\alpha})$. This boils down to checking \eqref{hopfcocycle} for $\tau.\alpha$ with the $\tl^{\alpha}$-action: 
\begin{align*}(\tau.\alpha)(xy)&=\tau(\alpha(x)\alpha(y))=\tau(\alpha(x))(\alpha(x)\bullet\tau(\alpha(y))\\&=\tau(\alpha(x))(x\bullet^{\alpha}\tau(\alpha(y))=(\tau.\alpha)(x)(x\bullet^{\alpha}(\tau.\alpha)(y)).\end{align*}
As for $\mathbb{A}(\tl)$-linearity, for every $\phi\in\mathbb{A}(\tl)$, we have
\begin{align*}((\tau.\alpha).\phi)(x,a,b)&=(\tau.\alpha)(x,a\phi^{-1},b\phi^{-1})=\tau(\alpha(x),a\phi^{-1},b\phi^{-1})\\&=(\tau.\phi)(\alpha(x),a,b)=((\tau.\phi).\alpha)(x,a,b).\end{align*}\qed

We need several short remarks.
\begin{Lem}\label{correction} Suppose $\tau$ is a $2$-cocycle. Assume $r\in\bu{(\D{G})}$ is such that $\phi_p.r=\epsilon$. Set $r_i=\phi_i.r,\,1\le i\le p$. Define a $1$-cocycle $\zeta:C_p\to\bu{(\k^G)}$ by $\zeta(t^i)=r_i$ and a $2$-cocycle $\tau'=\tau(\delta_{G}\zeta)$. Then the mapping
\begin{equation*}\iota:H(\tau,\tl)\to H(\tau',\tl),\,\iota(p_at^i)=p_ar_it^i,\,a\in G,\,1\le i\le p\end{equation*}
is an equivalence of extensions.\end{Lem}
\pf It suffices to show show $\delta_{G}\zeta\in B^2_c$ for then \cite[5.2]{M} yields the conclusion of the lemma. Now  $\delta_{G}\zeta\in B^2_c$ means that $\zeta$ satisfies \eqref{hopf1cocycle}. The argument of Lemma \ref{hf=adm} used to derive \eqref{hopfcocycle} from the condition \eqref{Adm} works verbatim for $\zeta$.\qed

\begin{Lem}\label{cocommutative} $H(\tau,\tl)$ is cocommutative  iff $\tau$ lies in $H^2_{\c}(\tl)$.\end{Lem}
\pf $H^*(\tau,\tl)$ is commutative iff $\ov{a}\ov{b}=\ov{b}\ov{a}$ which is equivalent to $\tau(a,b)=\tau(b,a)$. This condition is equivalent to  $\tau(t):G\times G\to\bu{\k}$ being a symmetric $2$-cocycle. Indeed, one implication is trivial, while if $\tau(t)$ is symmetric, then as pointed out in the odd case of Proposition \ref{hopf2cohomology} $\tau(t)$ is a coboundary, that is an element of $B^2_N/\ker\Phi$. A reference to Lemma \ref{structureofcoboundaries}(i) completes the proof.\qed

Unless stated otherwise, $H(\tau,\tl)$ is a noncocommutative Hopf algebra. We pick another algebra $H(\tau',\tl')$ isomorphic to $H(\tau,\tl)$ via $\psi:H(\tau,\tl)\to H(\tau',\tl')$. The next observation is noted in \cite[p. 802]{Mas1}. 
\begin{Lem}\label{Gstability} Mapping $\psi$ induces an Hopf automorphism of $\D{G}$.\end{Lem}\qed

Let $G$ be a finite group and $\mathrm{Aut_{Hf}}(\D{G})$ be the group of Hopf automorphisms of $\D{G}$. Identifying $(\k^G)^*$ with $\k G$ as in \S\ref{abelian}, for every $\phi\in\mathrm{Aut_{Hf}}(\D{G})$ the transpose mapping $\phi^*$ is a Hopf automorphism of $\k G$, hence an automorphism of $G$. This leads up to
\begin{Lem}\label{dualofphi} Let $G$ be a finite group. The mapping $\phi\mapsto\phi^*$ is an isomorphism between $\mathrm{Aut_{Hf}}(\D{G})$ and $\mathrm{Aut}(G)$. $\phi$ is a $C_p$-isomorphism $(\D{G},\bullet)\to (\D{G},\circ)$ if and only if $\phi^*$ is a $C_p$-isomorphism $(G,\tl')\to (G,\tl)$.\end{Lem}
\pf The first assertion is clear by the opening remark. Next we recall that $\phi^*$ acts on $G$ via
\begin{equation}\label{dualphi}(a\phi^*)(f):=f(a\phi^*)=\phi(f)(a),\,f\in\D{G}.\end{equation} 
Let $\tl$ and $\tl'$ be actions related to $\bullet$ and $\circ$ by \eqref{dualFaction}. The last conclusion follows from the calculation 
\begin{align*} ((a\tl' x)\phi^*)(f)&=\phi(f)(a\tl' x)=(x\circ\phi(f))(a)=\phi(x\bullet f)(a)\\
&=(a\phi^*)(x\bullet f)=(a\phi^*\tl x)(f),\,\text{for all}\,f\in\D{G}.\end{align*}\qed

We proceed to the formulation of isomorphism theorems. First we rephrase definitions of $[\tl]$ and $C(\tl)$. Let `$\simeq$' denote equivalence of actions of $C_p$ on $G$. With $\mathcal{R}$ defined in the Introduction we have $[\tl]=\{\tl'\in\mathcal{R}|\tl'\simeq\tl^{\alpha}\,\text{for some}\,\alpha\in A_p\}$ and $C(\tl)=\{\alpha\in A_p|\tl^{\alpha}\simeq\tl\}$. Furthermore we denote by $\gg$ the subgroup of $\mathrm{Aut}(G)$ generated by $\A$ and a set of automorphisms $\lambda_{\alpha}\in I(\tl,\tl^{\alpha})$ one for every $\alpha\in C(\tl)$ if $\tl$ is nontrivial, and $\A\times A_p$, otherwise. 
\begin{Prop}\label{keygroup} $\g(\tl)$ is a crossed product of $\A$ with $C(\tl)$.\end{Prop}
\pf The claim holds by definition for the trivial action. Else, we note that $\lambda\A\lambda^{-1}=\A$ for every $\lambda\in I(\tl,\tl^{\alpha})$ by Lemmas \ref{actiononcocycles}(ii),\ \ref{secondaction}(i). In addition, for every $\lambda,\mu\in I(\tl,\tl^{\alpha}),\,\lambda^{-1}\mu\in\A$. Thus we have $I(\tl,\tl^{\alpha})=\A\lambda_{\alpha}$. It follows that $\lambda_{\alpha}\cdot\lambda_{\beta}=\phi(\alpha,\beta)\lambda_{\alpha\beta}$ for some $\phi(\alpha,\beta)\in\A$. It remains to show that the kernel of $\pi:\g(\tl)\to C(\tl),\pi(\phi\lambda_{\alpha})=\alpha$ equals $\A$. Pick $\alpha:x\to x^k,k\ne 1$. Clearly $\lambda\in I(\tl,\tl^{\alpha})$ iff $t\lambda=\lambda t^k$ where we treat $t\in C_p$ as automorphism of $G$. Since elements of $\A$ commute with $t$, $I(\tl,\tl^{\alpha})\cap\A=\emptyset$.\qed

Our next goal is to define a $\g(\tl)$-module structure on $H^2_c(\tl)$. As we mentioned above $H^2_c(\tl)$ is $\A$-module. Further, for every $\lambda\in I(\tl,\tl^{\alpha})$ Lemmas \ref{actiononcocycles}(ii),\ \ref{secondaction}(ii) show that the mapping
\begin{equation}\label{mainautomorphisms}\omega_{\lambda,\alpha}: \tau\mapsto\tau.\lambda\alpha^{-1},\tau\in H^2_c(\tl)\end{equation}
is an automorphism of $H^2_c(\tl)$. We denote by $\ov{\phi}$ the automorphism of $H^2_c(\tl)$ induced by $\phi\in\A$ and we abbreviate $\omega_{\lambda_{\alpha},\alpha}    $ to $\omega_{\alpha}$.
\begin{Lem}\label{gactiononcohomology} The mapping $\phi\lambda_{\alpha}\mapsto\ov{\phi}\omega_{\alpha},\phi\in\A,\alpha\in C(\tl)$ defines $\g(\tl)$-module structure on $H^2_c(\tl)$.\end{Lem}
\pf $H^2_c(\tl)$ is a subquotient of $Z^2(G,\bu{(\k^{C_p})})$, and the action of $\A$ and $\omega_{\alpha}$ on $H^2_c(\tl)$ are induced from their action on $Z^2(G,\bu{(\k^{C_p})})$. Furthermore it is elementary to check that every $\lambda\in\mathrm{Aut}(G)$ commutes with every $\beta\in A_p$ as mappings of $Z^2(G,(\k^{C_p})^{\bullet})$. It follows that the equalities $\omega_{\alpha}\omega_{\beta}=\ov{\phi(\alpha,\beta)}\omega_{\alpha\beta}$ and $\omega_{\alpha}\ov{\phi}\omega_{\alpha}^{-1}=\lambda_{\alpha}\ov{\phi}\lambda_{\alpha}^{-1}$ hold in $\mathrm{Aut}(Z^2(G,\bu{(\k^{C_p})})$. This shows that the mapping of the Lemma is a homomorphism, as needed.\qed
\begin{Thm}\label{mainthm} {\rm (I)}. Noncocommutative extensions $H(\tau,\tl)$\\ and $H(\tau',\tl')$ are isomorphic if and only if 
\begin{enumerate}
\item[(i)] There exist an $\alpha\in A_p$ and a $C_p$-isomorphism $\lambda:(G,\tl)\to(G,\tl'^{\alpha})$ such that 
\item[(ii)] $\tau'=\tau.(\lambda\alpha^{-1})$ in $H^2_c(\tl')$.
\end{enumerate}
{\rm(II)}. There is a bijection between the orbits of $\g(\tl)$ in $H^2_c(\tl)$ not contained in $H^2_{\c}(\tl)$ and the isomorphism types of noncocommutative extensions in $\mathrm{Ext}_{[\tl]}(\k C_p,\k^G)$.
\end{Thm}
\pf (I). In one direction, suppose $\psi:H(\tau,\tl)\to H(\tau',\tl')$ is an isomorphism. By Lemma \ref{Gstability} $\psi$ induces  an automorphism $\phi:\D{G}\to\D{G}$, and from  Lemma \ref{gradedpsi} we have the equality $\psi(t)=rt^k$ for some $k$ and $r\in\D{G}$. The equality $\psi(t^p)=1$ implies $(rt^k)^p=\phi_p(t^k)\circ r=1$ and, as $\phi_p(t^k)=\phi_p(t)$, we have $\phi_p\circ r=1$ which shows $r\in\bu{(\D{G})}$. Let $\alpha:x\mapsto x^k, x\in C_p$ be this automorphism of $C_p$, and set $\phi=\psi|_{\k^G}$. Then the calculation
\begin{equation*}\phi(t\bullet f)=\psi(tft^{-1})=r\alpha(t)\phi(f)\alpha(t)^{-1}r^{-1}=\alpha(t)\circ\phi(f),\,f\in\D{G}\end{equation*}
shows $\phi:(\D{G},\bullet)\to(\D{G},\circ^{\alpha})$ is a $C_p$-isomorphism. It follows by Lemma \ref{dualofphi} that $(G,{\tl'}^{\alpha})$ is isomorphic to $(G,\tl)$ under $\phi^*$, hence $\lambda=(\phi^*)^{-1}:(G,\tl)\to (G,{\tl'}^{\alpha})$ is a required isomorhism.

It remains to establish the second condition of the theorem. To this end we first modify $\psi$. Namely, set $s=\phi^{-1}(r)$ and observe that, as $\phi^{-1}$ is a $C_p$-mapping and $\phi_p\circ^{\alpha}r=1$, we get $\phi_p\bullet s=\phi^{-1}(\phi_p\circ^{\alpha}r)=1$. Therefore by Lemma \ref{correction} there is an equivalence $\iota: H(\tau,\tl)\to H(\widetilde{\tau},\tl)$ with $\iota(t)=st$. Notice that $\iota$ is an algebra map with $\iota(s)=s$ for all $s\in\k^G$, hence $\iota^{-1}(t)=s^{-1}t$. Thus we have $(\psi\iota^{-1})(t)=t^k$ by the choice of $s$. It follows we can assume $\psi(t)=t^k$ hence $\psi(x)=x^k$ for all $x\in C_p$. 

Abbreviating $H(\tau,\tl),H(\tau',\tl')$ to $H,H'$, respectively, we take up the identity. 
\begin{equation*}\Delta_{H'}(\psi(x))=(\psi\o\psi)\Delta_H(x), x\in C_p,\end{equation*}
expressing comultiplicativity of $\psi$ on elements of $C_p$. By \eqref{DeltaH} this translates into
\begin{equation}\label{identity}\sum_{a,b}\tau'(x^k,a,b)p_ax^k\o p_bx^k=\sum_{c,d}\tau(x,c,d)\phi(p_c)x^k\o\phi(p_d)x^k.\end{equation}
Next we connect $\phi(p_b)$ to the action of $\phi^*$. This is given by the formula 
\begin{equation}\label{actionofphi}\phi(p_b)=p_{b(\phi^*)^{-1}}.\end{equation}

For, since $\phi$ is an algebra map, $\phi(p_b)=p_c$ where $c$ is such that $\phi(p_b)(c)=1$. By definition of action $\phi^*$, $\phi(p_b)(c)=(c\phi^*)(p_b)=p_b(c\phi^*)$, hence $c\phi^*=b$, whence $c=b(\phi^*)^{-1}$.

Switching summation symbols $c,d$ to $l=c(\phi^*)^{-1}$ and $m=d(\phi^*)^{-1}$, the right-hand side of \eqref{identity} takes on the form
\begin{equation*} \sum_{l,m}\tau(x,l\phi^*,m\phi^*)p_lx^k\o p_mx^k\end{equation*}
Thus $\psi$ is comultiplicative on $C_p$ iff 
\begin{equation}\label{comultiplicative}\tau'(\alpha(x),a,b)=\tau(x,a\phi^*,b\phi^*)=\tau(\phi^*)^{-1}(x,a,b)=\tau.\lambda(x,a,b).\end{equation}
Applying $\alpha^{-1}$ to the last displayed equation we arrive at
\begin{equation}\label{comultiplicative'}\tau'(x,a,b)=\tau.\lambda\alpha^{-1}(x,a,b).\end{equation}
as needed. 

Conversely, let us assume hypotheses of part (I). Using Lemma \ref{dualofphi} we infer that $\lambda^{-1}$ induces a Hopf $C_p$-isomorphism

\noindent $\phi=(\lambda^{-1})^*:(\D{G},\bullet)\to(\D{G},\circ^{\alpha})$. We define 
\begin{equation*}\psi:H(\tau,\tl)\to H(\tau',\tl')\,\text{via}\,\psi(fx)=\phi(f)\alpha(x),f\in\D{G},x\in C_p.\end{equation*}
First we verify that $\psi$ is an algebra map utilizing $\phi(x\bullet f)=\alpha(x)\circ\phi(f)$, namely 
\begin{align*}&\psi((fx)(f'x'))=\psi(f(x\bullet f')xx')=\phi(f)\phi(x\bullet f')\alpha(x)\alpha(x')\\
&=\phi(f)(\alpha(x)\circ\phi(f'))\alpha(x)\alpha(x')=\phi(f)\alpha(x)\phi(f')\alpha^{-1}(x)\alpha(x)\alpha(x')\\
&=(\phi(f)\alpha(x))(\phi(f')\alpha(x')=\psi(fx)\psi(f'x').\end{align*}
To see comultiplicativity of $\psi$ we need to verify 
\begin{equation}\label{comult}\Delta_{H'}(\psi(fx))=(\psi\otimes\psi)\Delta_H(fx).\end{equation}
By the multiplicativity of $\Delta_{H'},\psi,\Delta_H$ it suffices to check \eqref{comult} separately for any $f$ and for every $x$. Now the first case holds as $\phi$ is a coalgebra mappping, and the second follows from $\tau'=\tau.\lambda\alpha^{-1}$ by calculations \eqref{identity} and \eqref{comultiplicative'}.

(II). Let $\mathcal{H}$ denote the set of all pairs $(\tau',\tl')$ with $\tl'$ and $\tau'$ running over $[\tl]$ and $H^2_c(\tl')\setminus H^2_{\c}(\tl')$, respectively. We define an equivalence relation on $\mathcal{H}$ by 
\begin{equation*}(\tau',\tl')\sim(\tau'',\tl'')\;\text{iff}\; H(\tau',\tl')\simeq H(\tau'',\tl'').\end{equation*}
Let $\mathcal{H}/\sim$ stand for the set of equivalence classes. By construction $\mathcal{H}/\sim$ is just a copy of ${}_{\text{nc}}\mathrm{Ext}_{[\tl]}(\k C_p,\k^G)/\cong$. We select the subset $\mathcal{H}(\tl)=\{(\tau,\tl)|\tau\in H^2_c(\tl)\}$ of $\mathcal{H}$ and define the orbit $(\tau,\tl)\gg$ as the set $\{(\tau',\tl)|\tau'\in\tau\gg\}$. The proof will be complete if we show that the set $\{C\cap\mathcal{H}(\tl)|C\in\mathcal{H}/\sim\}$ coincides with the set of orbits of $\gg$ in $\mathcal{H}(\tl)$. Now pick $(\tau',\tl')\in C$. Since $\tl'\in[\tl]$, there exists an isomorphism $\mu\in I(\tl',\tl^{\alpha})$ hence setting $\tau=\tau'.\mu\alpha^{-1}$ we have $(\tau,\tl)\in C$ by part (I). Moreover $C\cap\mathcal{H}(\tl)\ni(\sigma,\tl)$ if and only if $H(\tau,\tl)\simeq H(\sigma,\tl)$, hence by part (I) again we have $\sigma=\tau.\omega_{\lambda,\alpha}$, that is $\sigma\in\tau\gg$. Same argument shows that the equivalence class generated by $(\tau,\tl)$ intersect $\mathcal{H}(\tl)$ in the orbit of $(\tau,\tl)$.\qed

\begin{Cor}\label{sizeoforbit} For every $\tau\in H^2_c(\tl)$ the cardinality of the orbit $\tau\gg$ satisfies
\begin{equation*}|\tau\A|\le|\tau\g(\tl)|\le|C(\tl)||\tau\A|.\end{equation*}
\end{Cor}
\pf Since $\A\subseteq\g(\tl)$ the lower bound is clear. By Propositon \ref{keygroup} $\g(\tl)=\bigcup_{\alpha\in C(\tl)}\omega_{\alpha}\A$ hence $\tau\g(\tl)=\bigcup_{\alpha\in C(\tl)}\tau\omega_{\alpha}\A$. It remains to note that for every $\alpha$ the cardinality of $\tau\omega_{\alpha}\A$ coincides with that of $\tau\A$.\qed

With some extra effort we can extend the bijection theorem to the entire set $\E$ provided $G$ is an elementary $p$-group for any $p$. Since our prime interest lies with nontrivial Hopf algebras we state the result without proof. 
\begin{Thm}\label{maincocom} Let $G$ be a finite elementary $p$-group. The number of isotypes of cocommutative Hopf algebras in $\E$ equals to the number of orbits of $\A$ in  $H^2_{\c}(\tl)$.\end{Thm}\qed

We comment briefly on the dual case of commutative Hopf algebras. First, $\E$ contains a commutative Hopf algebra iff $\tl=\text{triv}$. Second, we introduce the group $\mathrm{Cext}(G,C_p)$ of central extensions of $G$ by $C_p$ \cite{BT}. We outline properties of $\mathrm{Ext}_{[\t]}(\k C_p,\k^{G})$ again without proof.
\begin{Thm}\label{maincom} {\rm(1)} The group $\mathrm{Ext}_{[\text{\rm triv}]}(\k C_p,\k^{G})$ is isomorphic to the group $\mathrm{Cext}(G,C_p)$  under the map $H(\tau,\text{\rm triv})\leftrightarrows \k^{L(\tau)}$ where $L(\tau)$ is the central extension defined by the $2$-cocycle $\tau$.

{\rm(2)} For $G$ elementary $p$-group of rank $n$ with an odd $p$ the number of isotypes in $\mathrm{Ext}_{[\text{\rm triv}]}(\k C_p,\k^G)$ equals $\lfloor\frac{3n+2}{2}\rfloor$.
\end{Thm}

For calculation of orbits of $\g(\tl)$ in $H^2_c(\tl)$ we prefer to use its isomorphic copy $Z^2_N(\tl)/\ker\Phi$ which we denote by $\X$ and refer to it as the classifying group for $\E$.

We turn $\mathbb{X}(\tl)$ into a $\g(\tl)$-module by transfering its action from $Z^2_c(\tl)$ to $Z^2_N(\tl)$ along $\Theta$. Pick some $\omega_{\lambda,\alpha}$ and suppose $\alpha^{-1}:x\mapsto x^l,x\in C_p$. For $s\in Z^2_N(\tl)$ we put
\begin{equation}\label{derivedaction}s.\omega_{\lambda,\alpha}=(\phi_l\bullet s).\lambda.\end{equation} 
\begin{Lem}\label{gisomorphism} {\rm(i)} For every prime and any action `$\tl$' the isomorphism $\Theta_*:H^2_c(\tl)\simeq\X$ of Corollary \ref{basicisic} is  $\g(\tl)$-linear.\\ 
{\rm(ii)} For every prime and any action $\X$ fits into the exact sequence
\begin{equation}\label{seqforX}\G^{C_p}/N(\G)\rightarrowtail\X\twoheadrightarrow\un{a}(Z^2_N(\tl)).\end{equation} 
{\rm(iii)} For every odd $p$ there is a $\gg$ splitting
\begin{equation}\label{Xsplitting} \X\simeq \G^{C_p}/N(\G)\times\mathrm{Alt}_N(G).\end{equation}
\end{Lem} 
\pf (i) We begin by noting that for every $\lambda\in I(\tl,\tl^{\alpha})$ there holds (*) $x\bullet^{\alpha}(s.\lambda)=(x\bullet s).\lambda,\, x\in C_p$. Still assuming $\alpha^{-1}:x\to x^l$, the conclusion (i) follows from \eqref{componentsoftau} and the opening remark by the calculation
\begin{align*}&\Theta(\tau.\omega_{\lambda,\alpha})=(\tau.\omega_{\lambda,\alpha})(t)=(\tau.\lambda)(t^l)=(\text{by}\,\eqref{componentsoftau})\,\phi_l\bullet^{\alpha}(\tau.\lambda)(t)\\&=\phi_l\bullet^{\alpha}(\tau(t).\lambda)=(\text{by (*)})\,(\phi_l\bullet\tau(t)).\lambda=\Theta(\tau).\omega_{\lambda,\alpha}\end{align*}
This equation demonstrates that definition \eqref{derivedaction} turnes $Z^2_N(\tl)$ into a $\g(\tl)$-module. It is immediate that $B^2_c(\tl)$ is a $\g(\tl)$-subgroup of $Z^2_c(\tl)$. By Lemma \ref{structureofcoboundaries}(ii) $\ker\Phi$ is a $\g(\tl)$-subgroup, which proves part (i).

(ii) The mapping $\un{a}:Z^2(G,\k^{\b})\to\mathrm{Alt}(G)$ of Proposition \ref{hopf2cohomology} restricted to $Z^2_N(\tl)$ gives rise to an exact sequence $B^2_N(\tl)\to Z^2_N(\tl)\to\un{a}(Z^2_N(\tl)$. Thanks to the $\gg$-isomorphism $\G^{C_p}/N(\G)\simeq B_N^2(\tl)/\ker\F$ induced by $\F$ (see Lemma \ref{structureofcoboundaries}) we arrive at the exact sequence \eqref{seqforX} of $\gg$-modules. 

(iii) For an odd $p$ splitting \eqref{Z^2adm} is carried out by the mapping $s\mapsto s\un{a}(s^{-2})\times\un{a}(s^2)$ which is clearly a $\g(\tl)$-map. It remains to note that homomorphism $\Phi$ is also a $\g(\tl)$-map.\qed

We point out that part (ii) fails in general for $2$-groups.\footnotemark\footnotetext{ See Appendix 2} 


\section{Almost Abelian Hopf Algebras of Dimension $\le p^4$}\label{p^4}

\subsection{Hopf algebras of dimension $\le p^3$} We begin by revisiting classification of semisimple Hopf algebras of dimension $p^2,p^3$ due to \cite{{Mas3},{Mas1}}. If $\dim H=p^2$, then by a Kac-Masuoka theorem \cite{{Kac},{Mas3}} $H$ contains a central subHopf algebra $\k C_p$ hence $H\in\mathrm{Ext}_{[\t]}(\k C_p,\k^{C_p})$. Thus $H$ is commutative, and as $\text{Alt}(C_p)=1$, $H$ is cocommutative. It follows that $H=\k L$ where $L$ is a group of order $p^2$, that is $L=C_{p^2}$ or $C_p\times C_p$.  

Suppose $\dim H=p^3$. By the Kac-Masuoka theorem, loc.cit, applied to $H^*$ we have that $H^*$ is a central extension of the form $\k C_p\rightarrowtail H^*\twoheadrightarrow Q$ where $\dim Q=p^2$. By the foregoing $Q=\k G$ with $G=C_{p^2}$ or $C_p\times C_p$. By duality $H$ is a cocentral extension of $\k C_p$ by $\k^G$. If $G=C_{p^2}$, then $\text{Alt}(G)=1$, hence $H$ is cocommutative. It follows that a nontrivial $H$ belongs to $\mathrm{Ext}(\k C_p,\k^{C_p\times C_p})$ with a nontrivial action of $C_p$ on $C_p\times C_p$.

Before moving on we introduce algebras $R_i=\Z_pC_p/\gen{(t-1)^i},\,0\le i\le p-1 $ and make a notational change. Below we write $\alpha_k$ for the mapping $x\mapsto x^k,x\in C_p$, $\tl^k$ for $\tl^{\alpha_k}$ and $\omega_k$ for $\omega_{\alpha_k}$. The arguments in the next proposition will be used throughout \S 4.2.2.
\begin{Prop}{\rm(\cite{Mas1})}\label{C_ptimesC_p} There are up to isomorphism $p+7$ Hopf algebras in $\mathrm{Ext}(\k C_p,\k^{C_p\times C_p})$, $p+1$ of which are nontrivial.\end{Prop}
\pf  We run the procedure for computing the number of isoclasses for $G=C_p\times C_p$. Let $\tl_r$ denote the right regular action of $C_p$ on $R_2$. Every nontrivial $C_p$-module $(C_p\times C_p,\tl)$ is isomorphic to $(R_2,\tl_r)$. In consequence 
$\mathrm{Ext}(\k C_p,\k^G)=\mathrm{Ext}_{[\tl_r]}(\k C_p,\k^G)\cup \mathrm{Ext}_{[\t]}(\k C_p,\k^G)$. By Theorem \ref{maincom}(2) $\mathrm{Ext}_{[\t]}(\k C_p,\k^G)$ contributes four nonisomorphic algebras. It remains to show that $\mathrm{Ext}_{[\tl_r]}(\k C_p,\k^G)$ contains $p+3$ isotypes. To simplify notation we put $\tl=\tl_r$.

(i) The classifying group $\mathbb{X}(\tl)$.  Set $G=R_2$ and let $e=\ov{1},f=\ov{t-1}$ where $\ov{r}$ is the image of $r\in R_0$ in $R_2$. The matrix of $t$ in the basis $\{e,f\}$ is $T=\begin{pmatrix}1&1\\0&1\end{pmatrix}$. Let $\{e^*,f^*\}$ be the dual basis for $\G$. The mapping induced by $t$ in $\G$ has the matrix $T^{\Tr}$ relative to the dual basis. Hence $e^*$ is fixed by $t$ and $N(\G)=(t-1)^{p-1}.\G=0$, as $p>2$, in the additive notation. It follows that $\G^{C_p}/N(\G)=\gen{e^*}$. Further $\mathrm{Alt}(G)=\G\we\G$, where $\wedge$ denote the multiplication in the Grassman algebra over $\G$, is generated by $e^*\we f^*$, and the latter is a fixed by $t$. Therefore $\phi_p(t).e^*\we f^*=p(e^*\we f^*)=0$, hence $\mathrm{Alt}_N(G)=\mathrm{Alt}(G)$. All in all we arrive at the equality
\begin{equation*}\mathbb{X}(\tl)=\gen{e^*,e^*\we f^*}\end{equation*} 

(ii) Groups $\A,C(\tl)$ and $\gg$. By definition $\phi\in\A$ iff the matrix $\F$ of $\phi$ satisfies $\F T=T\F$ and $\det\F\ne 0$. This condition is equivalent to $\F=\begin{pmatrix}c&d\\0&c\end{pmatrix}, c\in\mathbb{Z}_p^{\bullet}$. By the opening remark $C(\tl)=\mathbb{Z}_p^{\bullet}$ as $(G,\tl^k)\simeq(G,\tl)$ for every $k\in\mathbb{Z}_p^{\bullet}$.

The group $\gg$ is generated by $\A$ and a set $\{\lambda_k|k\in\mathbb{Z}_p^{\bullet}\}$ with $\lambda_k\in I(\tl,\tl^k)$. An easy verification gives that $\lambda_k$ 
defined via 

\noindent$e.\lambda_k=e,f.\lambda_k=kf$ lies in $I(\tl,\tl^k)$.

(iii) Orbits of $\gg$ in $\X$. First we determine the orbits of $\A$. Pick $\phi\in\A$ and suppose it has the matrix $\F$ relative to $\{e,f\}$. It is an elementary fact that the mapping induced by $\phi$ in $\G$ has the matrix $\F^{\text{tr}}$ in the dual basis. If $\F$ is written as in (ii) then we have
\begin{equation*} e^*.\phi^{-1}=ce^*,\;\text{and}\;e^*\we f^*.\phi^{-1}=c^2e^*\we f^*.\end{equation*}
Let us identify $ae^*+be^*\we f^*\in\X$ with the vector $(a,b)\in\mathbb{Z}_p^2$. By the above $\phi\in\A$ acts in $\mathbb{Z}_p^2$ via $(a,b).\phi^{-1}=(ca,c^2b)$.

By Corollary \ref{sizeoforbit} the $\gg$-orbit of $(a,b)$ is the union of $\A$-orbits of elements $(a,b).\omega_k,k\in C(\tl)$. There $\omega_k=\lambda_k\alpha_k^{-1}$, and for every $x\in\X$ there holds by \eqref{derivedaction} $x.\omega_k=(\phi_l.x).\lambda_k$ where $l=k^{-1}$. Since $e^*$ and $e^*\we f^*$ are fixed by $t$ we have $\phi_l.x=lx$ for $x=e^*,e^*\we f^*$. Moreover it is immediate that $e.\lambda_k=e^*$ and $e^*\we f^*.\lambda_k=le^*\we f^*$. We conclude that $(a,b).\omega_k=(la,l^2b)\in(a,b)\A$. It follows that $\gg$-orbits coincide with $\A$-orbits. We compute the latter.

The subset ${\Z_p^{\bullet}}^2$ of $\Z_p^2$ is stable under action of $\A$. For every $m\in\Z_p^{\bullet}$ the set $(1,m)\A$ has $p-1$ elements and moreover $(1,m)\A\cap(1,n)\A=\emptyset$ if $m\ne n$. Since $|{\Z_p^{\bullet}}^2|=(p-1)^2$ the family 

\noindent$\{(1,m)\A|m\in\Z_p^{\bullet}\}$ accounts for all orbits in ${\Z_p^{\bullet}}^2$. Thus we obtained $p-1$ nontrivial orbits. The complement $\Z_p^2\setminus{\Z_p^{\bullet}}^2$ is the union of $\{(0,b)|b\in\Z_p^{\bullet}\}$ and $\{(a,0)|a\in\Z_p\}$. Let $\zeta$ be a generator of $\Z_p^{\bullet}$. It follows readily that $\{(0,b)|b\in\Z_p^{\bullet}\}$ is the union of $(0,1)\A$ and $(0,\zeta)\A$ which supplies two more nontrivial orbits. The second set is the union of two trivial orbits, viz. $\{(0,0)\}$ and its complement. \qed

To recover the full strength of \cite{Mas1} we would need to show that every $H(\tau,\tl)$ is self-dual. However, such a theorem is unattainable due to the next
\begin{Rem}\label{counterexample} Let $\tau(t)=e^*\we f^*$ and $H(\tau,\tl)$ be the corresponding Hopf algebra. $H(\tau,\tl)^*\simeq H(\tau,\tl)$ if and only if $\frac{p-1}{2}$ is a square in $\mathbb{Z}_p^{\b}$.\end{Rem} 
\pf By general theory $H(\tau,\tl)^*\simeq H(\tau',\tl)$ for some $\tau'\in H^2_c(\tl)$. An isomorphism $H(\tau,\tl)^*\simeq H(\tau,\tl)$ exists if and only if $\tau'(t)$ and $\tau(t)$ lie on the same orbit. The $2$-cocycle $\tau'$ is the multiplication cocycle for $H(\tau,\tl)$ written as an element of $\mathrm{Ext}(\k G,\k^{C_p})$. Since $H(\tau,\tl)=\k(\G\rtimes C_p)$ and $e^*$ is a fixed point under the action of $C_p$, $\k\gen{e^*}$ is a normal subHopf algebra of $H(\tau,\tl)$ giving rise to an exact sequence
\begin{equation}\label{seq}\k\gen{e^*}\hookrightarrow H(\tau,\tl)\overset\Pi\twoheadrightarrow \k\ov{G}\end{equation}
where $\ov{G}=\gen{x,y}$ with $x=\Pi(f^*),y=\Pi(t)$. Clearly $xy=yx$ so that $\ov{G}=G$. Let $\rho_{\Pi}=(\text{id}\o\Pi)\Delta_H:H(\tau,\tl)\to H(\tau,\tl)\o\k\ov{G}$ be the coaction induced by $\Pi$. We want to find a section $\gamma:\k\ov{G}\to H(\tau,\tl)$ splitting \eqref{seq}. This is the matter of finding $T$ satisfying $\rho_{\Pi}(T)=T\o y$. It is not hard to see that $T$ must be of the form $ut$ for some unit $u\in\k^G$in fact a tedious but straightforward verification shows that for $u=\sum_{i,j}\zeta^{-ij}p_{e^if^j}$\  $T=ut$ is a desired element. Since $f^*$ is a group-like element of $H(\tau,\tl)$, $\rho_{\Pi}({f^*}^iT^j)={f^*}^iT^j\o x^iy^j$, and therefore $\gamma:x^iy^j\mapsto {f^*}^iT^j$ defines a section of $\k\ov{G}$ in $H(\tau,\tl)$. Let $\tau':\ov{G}\times\ov{G}\to\k\gen{e^*}$ be the $2$-cocycle associated to $\gamma$. By definition $\tau'(a,b)=\gamma(a)\gamma(b)\gamma(ab)^{-1}$. We will write below $a=x^iy^j,b=x^ky^l$. A simple calculation using $Tf^*=f^*e^*T$ gives $\tau'(a,b)={e^*}^{jl}$. Viewing $e^*$ as the functional $e^*(t^k)=\zeta^k$ on $C_p(t)$ we conclude that $\tau'(t,a,b)=\zeta^{jl}$. 

We need to find a decomposition of $\tau'(t)$ according to \eqref{Z^2adm}, that is $\tau'(t)=\un{b}\cdot\lambda$ with $\un{b}\in B^2_N(\tl)$ and $\lambda\in\mathrm{Alt}_N(G)$. Set $\beta=\un{a}(\tau'(t))$ and note that by the definition of $\un{a}$,  $\beta(a,b)=\tau'(t,a,b)\tau'(t,b,a)^{-1}$ which gives $\beta(a,b)=\zeta^{jk-il}$. Observe that $\beta=\un{a}(\lambda)=\lambda^2$, hence $\lambda=\beta^{\frac{1}{2}}$ and therefore $\un{b}=\tau'\beta^{-\frac{1}{2}}$. It follows that $\un{b}(a,b)=(\zeta^{\frac{1}{2}})^{jk+il}$. Let us select $f:\ov{G}\to\k^{\b},\,f(x^iy^j)=(\zeta^{-\frac{1}{2}})^{ij}$. A straightforward calculation produces the equality $\delta f(a,b)=\un{b}$. We want to find the image of $\un{b}$ in $\X$ under $\F$, that is $\F(\un{b})=\phi_p(t).t$. Since $\phi_p(t).f(a)=\prod_{i=0}^{{p-1}} f(a^i\tl t^i)$, $\phi_p(t).f(x)=1$ as $x\tl t=x$, and $\phi_p.f(y)=\prod_{i=0}^{{p-1}}f(a^ib)=\prod_{i=0}^{{p-1}}(\zeta^{-\frac{1}{2}})^i=1$. Since $\phi_p(t).f$ is a character of $\ov{G}$, $\phi_p(t).f=1$. Thus $\ov{b}\in\ker\F$ which means $\tau'(t)=\beta^{\frac{1}{2}}$ in $\X$. But $\beta=\tau^{-1}$ as $\tau(a,b)=(e^*\we f^*)(a,b)=\zeta^{il-jk}$. Thus $\tau'(t)=\tau(t)^{\frac{p-1}{2}}$ or $(\frac{p-1}{2})e^*\we f^*$ in the additive notation. By Proposition \ref{C_ptimesC_p}(iii) $\tau'(t)$ lies on the orbit of $\tau(t)$ iff $\frac{p-1}{2}$ is a square.\qed

\subsection{Hopf algebras of dimension $p^4$}. From now on we assume that $H$ is of dimension $p^4$ with an abelian group $G$ of grouplikes of order $p^3$. 
\begin{Thm}\label{mainp^4} There are $5p+23$ distinct nontrivial almost abelian Hopf algebras of dimension $p^4$ if $p>3$, $31$ if $p=3,e\ge 3$ and $33$, otherwise.\end{Thm}
\pf This will carried out in steps. In the additive notation $G=\Z_p^3$ or $G=\Z_{p^2}\oplus\Z_p$, and the theory splits into two parts.

\subsubsection{ $G=\Z_p^3$}.\label{CaseA}

There are up to isomorphism two nontrivial $\Z_pC_p$-module structures on $G$. Namely, if $C_p$-module $G$ is decomposable, then $G\simeq R_2\oplus R_1$, and $G\simeq R_3$, otherwise. 

(I) Suppose $G\simeq R_2\oplus R_1$, and let $\tl_d$ be the action of $C_p$ on $G$ composed of regular actions of $C_p$ on $R_2$ and $R_1$. We aim to prove
\begin{Thm}\label{isoclassesAI} $\mathrm{Ext}_{[\tl_d]}(\k C_p,\k^{C_p^3})$ contains $2p+11$ isotypes of extensions $2p+8$ of which are nontrivial.\end{Thm}
\pf We carry out the procedure for computing the number of isotypes for $C_p$-module $(G,\tl_d)$. To simplify notation we put $\tl=\tl_d$.
 
(1) The classifying group $\X$. Select a basis $\{e,g,f\}$ for $G$ where $\{e,f\}$ is the basis for $R_2$ as in Proposition \ref{C_ptimesC_p}, and $R_1=\Z_pg$. Clearly the matrix $T$ of $t$ in that basis is $T=\begin{pmatrix}1&0&1\\0&1&0\\0&0&1\end{pmatrix}$. Let $\{e^*,g^*,f^*\}$ be the dual basis for $\G$.  We fix a basis $\{e^*\we g^*,e^*\we f^*,g^*\we f^*\}$ for $\G\we\G$. We refer to the above bases as standard. 
\begin{Prop}\label{structureX1} $\X=\gen{e^*,g^*}\oplus\G\wedge\G$.\end{Prop} 
\pf Recall $\X= \G^{C_p}/N(\G)\oplus \text{Alt}_N(G)$. We use the well known identification  $\text{Alt}(G)=\G\we\G$. One can see easily that the matrix of $t$ in the standard basis of $\G$ is $T^{\Tr}$. By general principles \cite[III,8.5]{Bou} the matrix of $t$ in the standard basis of $\G\we\G$ is $T^{\Tr}\we T^{\Tr}=\begin{pmatrix}1&0&0\\0&1&0\\-1&0&1\end{pmatrix}$.
It follows that $(t-1)^{p-1}\bullet\G=0$ and $(t-1)^{p-1}\bullet\G\we\G=0$, that is $N(\G)=0$ and $(\G\we\G)_N=\G\we\G$. Further, one can see easily $\G^{C_p}=\gen{e^*,g^*}$.\qed

(2) Groups $\A,C(\tl)$ and $\gg$. By definition $\phi\in\A$ iff its matrix $\F$ satisfies $\F T=T\F$ and $\det\F\ne 0$. By a straighforward calculation one can see that $\phi\in\A$ iff  
\begin{equation}\label{A1}\F=\begin{pmatrix}a_{11}&a_{12}&a_{13}\\0&a_{22}&a_{23}\\0&0&a_{11}\end{pmatrix},\,a_{ij},\in\Z_p,\,a_{11}a_{22}\ne 0\end{equation}
It is easy to see that $(G,\tl^k)\simeq(G,\tl)$ for every $k\in\mathbb{Z}_p^{\bullet}$ which gives $C(\tl)=\bu{\Z_p}$. Likewise one can check directly that $\lambda_k:e\mapsto e, g\mapsto g,f\mapsto kf$ lies in $I(\tl,\tl^k)$ for every $k\in\bu{\Z_p}$. This determines $\gg$ as the latter is generated by $\A$ and the $\lambda_k$. 

(3) Orbits of $\A$ in $\X$. In order to simplify notation we change coordinates of matrices \eqref{A1} by setting $u=a_{11},v=a_{22}, a_{12}=u^{-1}q,\\a_{23}=u^{-1}r,a_{13}=s$. We treat the tuple $(u,v,q,r,s)$ as the coordinate of either $\phi$ or its matrix $\F$. On general principles \cite[III,8.5]{Bou} the matrices of $\phi^{-1}$ in the standard bases for $\G$ and $\G\we\G$ are $\F^{\Tr}$ and $\F^{\Tr}\we\F^{\Tr}$, respectively. For $\F=\F(u,v,q,r,s)$ a routine calculation gives
\begin{equation}\label{Finverse}\F^{\Tr}=\begin{pmatrix}u&0&0\\u^{-1}q&v&0\\s&u^{_1}r&u\end{pmatrix},\;\text{and}\end{equation} 
\begin{equation}\label{FwedgeF}\F^{\Tr}\we\F^{\Tr}=\begin{pmatrix}uv&0&0\\r&u^2&0\\z&q&uv\end{pmatrix},\end{equation}
where $z=\det\begin{pmatrix}u^{-1}q&v\\s&u^{-1}r\end{pmatrix}$. Next we identify $\X$ with $\mathbb{Z}_p^5$ via the assignment $x=a_1e^*+a_2g^*+b_1e^*\we g^*+b_2e^*\we f^*+b_3g^*\we f^*\mapsto v(x)=(a_1,a_2,b_1,b_2,b_3)$. We use the notation $e'_i,e''_j,i=1,2,j=1,2,3$ for the standard unit vectors in $\mathbb{Z}_p^2,\mathbb{Z}_p^3$, respectively. We begin with $\A$-orbits in $\G^{C_p}$ and $\G\we\G$. We define $Z'_i,Z''_j,\,0\le i\le 2,0\le j\le 3$ by the formula 
\begin{align*}Z'_i&=\{(a_1,a_2)|a_i\ne 0\,\text{and}\,a_k=0\,\text{for}\,k>i>0\},\\Z''_j&=\{(b_1,b_2,b_3)|b_j\ne 0\,\text{and}\,b_k=0\,\text{for}\,k>j>0\},\end{align*}
and $Z'_0=\{(0,0)\},Z''_0=\{(0,0,0)\}$. Furthermore we split $Z''_2$ into the union of $Z''_{2,k},k=0,1$ where $Z''_{2,k}=\{(b_1,\zeta^kb_2,0)|b_2\in{\bu{\Z_p}}^2\}$. We let $\q$ denote an element of $\{1,(2,0),(2,1),3\}$.

\begin{Lem}\label{partialorbits} The sets $Z'_i,Z''_{\q}$ are all the orbits of $\A$ in $\G^{C_p}$ and $\G\we\G$, respectively.\end{Lem}
\pf First note $\Z_p^2=\cup Z'_i$ and $\Z_p^3=\cup Z''_{\q}$. The equalities $e'_i\A=Z'_i,i=1,2$ are immediate by \eqref{Finverse}. This proves the first claim.

Similarly, using \eqref{FwedgeF} one can derive readily the equalities $e''_{\q}\A=Z''_{\q}$ for $\q=1,3$, and $\zeta^ke''_2\A=Z''_{2,k},k=0,1$.\qed

Let us write $Z'_i\times Z''_{\q}$ for the set of vectors $(v_1,v_2)$ with $v_1\in Z'_i,v_2\in Z''_{\q}$. These sets are $\A$-stable and some of them are orbits itself. We list those that are in
\begin{Lem}\label{splitorbits} For all $(i,\q)\ne(1,(2,k)),(2,3),k=0,1$ $Z'_i\times Z''_{\q}$ is an orbit.
\end{Lem}
\pf The claim is that for generators $e',e''$ of $Z'_i,Z''_{\q}$ in the nonexceptional cases, $(e',e'')$ generates $Z'_i\times Z''_{\q}$. We give details for $Z'_1\times    Z''_3$, other cases are treated similarly. Combining \eqref{Finverse} with \eqref{FwedgeF} we obtain 
\begin{equation*}(1,0,0,0,1).\A=\{(u,0,z,q,uv)\}\end{equation*}
Now for every element $(a_1,0,b_1,b_2,b_3)\in Z'_1\times Z''_3$ the equations $u=a_1, uv=b_3,uvr=b_1,q=b_2,$ are obviously solvable. A solution to the equation $z=b_1$ is provided by $r=0$ and $s=-v^{-1}b_1$. \qed

We pick up $p-1$ additional orbits in 
\begin{Lem}\label{extraorbits} Each set $Z'_1\times Z''_{2,k},k=0,1$ is a union of $(p-1)/2$ orbits.\end{Lem}
\pf Say $k=0$. By definition $Z'_1\times Z''_{2,0}=\{(a_1,0;b_1,b_2,0)|a_1\in\bu{\Z_p},b_2\in{\bu{\Z_p}}^2,b_1\,\text{arbitrary}\}$, hence $|Z'_1\times Z''_{2,0}|=\frac{(p-1)}{2}(p-1)p$. For every $m\in{\bu{\Z_p}}^2$ we let $z_m=(1,0;0,m,0)$. By \eqref{Finverse} and \eqref{FwedgeF} for $\phi=\phi(u,v,q,r,s)$  we have $z_m.\phi^{-1}=(u,0,mr,mu^2,0)$. A direct count gives $|z_m\A|=(p-1)p$, and one can verify directly that $z_m.\A\cap z_n.\A=\emptyset$ for $m\ne n$. Since there are $\frac{p-1}{2}$ orbits of this size, this case is done. For $i=1$ one should take $z'_m=(1,0;0,\zeta m,0)$.\qed

We summarize
\begin{Lem}\label{Aorbits}  There are $2p+8$ nontrivial orbits of $\A$ in $\X$.\end{Lem}
\pf The previous two lemmas give $p+8$ nontrivial orbits. The rest will come from splitting of the remaining set $Z'_2\times Z''_3$. The latter is defined as $\{(a_1,a_2,b_1,b_2,b_3)|a_2,b_3\in\bu{\Z_p},a_1,b_1,b_2\,\text{arbitrary}\}$. For every $k\in\Z_p$ we define $w_k=(k,1,0,0,1)$. Again by \eqref{Finverse} and \eqref{FwedgeF} we have
\begin{equation}\label{actionongenerators}  w_k.\phi^{-1}=(uk-u^{-1}q,v,z,q,uv).\end{equation}
where $(u,v,q,r,s)$ are the parameters of $\phi$. This formula shows that $w_k.\phi^{-1}$ does not depend on $r$, Setting $r=0$ we have $z=-sv$. It follows easily that $w_k.\phi^{-1}$ is uniquely determined by $(u,v,q,s)$, hence\\ $|w_k.\A|=(p-1)^2p^2$. Furthermore, we claim that $w_k.\A\cap w_l.\A=\emptyset$ for $k\ne l$. For, suppose 
\begin{equation*}(uk-u^{-1}q,v,-sv,q,uv)=(u'l-{u'}^{-1}q',v',-s'v',q',u'v')\end{equation*}
for some $(u,v,q,s)\,\text{and}\, (u',v',q',s')$. Then $v=v', q=q'$ give $u=u'$, hence $uk=ul$ and therefore $k=l$, a contradiction. We conclude that $|\cup_{0\le k\le p-1}w_k.\A|=p^3(p-1)^2$. As this is the number of elements in $Z'_2\times Z''_3$, the proof is complete.\qed 

(4) Orbits of $\gg$. We need to know the action of $\omega_k=\lambda_k\alpha_k^{-1}$ where $\lambda_k$ are defined in part (2). Set $l=k^{-1}\pmod p$.
\begin{Lem}\label{gaction} Action of $\omega_k$ is described by
\begin{align*} e^*.\omega_k&=le^*,\,g^*.\omega_k=lg^*\\
e^*\we g^*.\omega_k&=le^*\we g^*\\e^*\we f^*.\omega_k&=l^2e^*\we f^*\\g^*\we f^*.\omega_k&=-\binom l2 e^*\we g^*+l^2g^*\we f^*\end{align*}\end{Lem}
\pf By \eqref{derivedaction} for $x\in\X$, $x.\omega_k=(\phi_l\bullet x).\lambda_k$. For $x=e^*,g^*,e^*\we g^*,e^*\we f^*$ $\phi_l\bullet x=lx$ as these elements are fixed by $C_p$. Because $(t-1)^2\bullet\G\we\G=0$ we expand $\phi_l$ in powers of $t-1$, namely $\phi_l=l+\binom l2(t-1)+\text{higher terms}$. One can check $(t-1)\b g^*\we f^*=-e^*\we g^*$ which gives
\begin{equation*}\phi_l\bullet g^*\we f^*=lg^*\we f^*-\binom l2e^*\we g^*,\end{equation*}
By definition of $\lambda_k$ its matrix is $\Lambda_k=\text{diag}(1,1,k)$ (that is the diagonal matrix with entries $1,1,k$). It follows (see part (3)) that the matrix of $\lambda_k$ in the standard basis of $\G$ is $(\Lambda_k^{-1})^{\Tr}=\text{diag}(1,1,l)$. Applying $\lambda_k$ to $\phi_l\bullet x$ as $x$ runs over the standard bases of $\G$ and $\X$ we complete the proof of the Lemma.\qed

The next Proposition completes the proof of Theorem \ref{isoclassesAI}.
\begin{Prop}\label{GvsA} The sets of $\g(\tl)$ and $\A$-orbits coincide.\end{Prop}
\pf By Corollary \ref{sizeoforbit} for every $x\in\X$, $x\gg$ is a union of orbits $x.\omega_k\A$ for $1\le k\le p-1$. Thus it suffices to show  $x.\omega_k\in x\A$ for every $k$ and generators $x$ of every orbit of $\A$. We give a sample calculation for $x=w_m$ of Lemma \ref{Aorbits}. By Lemma \ref{gaction}
\begin{equation*} w_m.\omega_k=(lm,l,-\binom l2,0,l^2).\end{equation*} 
Now take $\phi$ with coordinates $u=l,v=l,q=0,r=0,s=l^{-1}\binom l2$. Then by \eqref{actionongenerators} $w_m.\phi^{-1}=w_m.\omega_k$ as needed.\qed. 

We move on to the next case

(II) $G\simeq R_3$. We denote by $\tl_r$ the right multiplication in $R_3$. This case is sensitive to the prime $p$. Let us agree to write $\x_p$ for $\mathbb{X}(\tl_r)$ if $G$ is a $p$-group. For $r\in\Z_pC_p$ we denote by $\ov{r}$ the image of $r$ in $R_3$. The elements $e=\ov{1},f=\ov{(t-1)},g=\ov{(t-1)^2}$ form a basis for $R_3$ in which action of $t$ is defined by $T=\begin{pmatrix}1&1&0\\0&1&1\\0&0&1\end{pmatrix}$. Let $\{e^*,f^*,g^*\}$ be the dual basis for $\G$, and $\{e^*\we f^*,e^*\we g^*,f^*\we g^*\}$ the induced basis for $\G\we\G$. We call all these bases standard. We aim to prove
\begin{Thm}\label{isoclassesAII}For $p>3$ $\mathrm{Ext}_{[\tl_r]}(\k^{C_p^3},\k C_p)$ contains $p+9$ isoclasses, $p+7$ of which are nontrivial, and three nontrivial isoclasses if $p=3$.\end{Thm}
\pf Proof will be carried out in steps following the procedure for computing the number of isoclasses. 

(1) Classifying groups $\x_p$. 
\begin{Lem}\label{structureX2} If $p=3$, then
\begin{equation*} \x_3=\gen{e^*\we f^*,e^*\we g^*}\end{equation*}
For every $p>3$
\begin{equation*} \x_p=\Z_pe^*\oplus\G\we\G\end{equation*}
\end{Lem}
\pf The matrices of $t$ in the standard bases of $\G$ and $\G\we\G$ are $T^{\Tr}$ and $T^{\Tr}\we T^{\Tr}$, respectively, with $T^{\Tr}\we T^{\Tr}=\begin{pmatrix}1&0&0\\1&1&0\\1&1&1\end{pmatrix}$. From this one computes directly $(t-1)^3\bullet\G=(t-1)^3\bullet\G\we\G=0$. Since $\phi_p(t)=(t-1)^{p-1}$, it follows that $N(G)=0$ and $(\G\we\G)_N=\G\we\G$ for any $p>3$. Furhermore $\G^{C_p}=\Z_pe^*$ for every $p$. Thus as $\x_p=\G^{C_p}/N(\G)\oplus(\G\we\G)_N$ the second statement of the Lemma follows.

Say $p=3$. Then $N(\G)=(t-1)^2\bullet\G=\Z_pe^*$, hence $\G^{C_p}/N(\G)=0$. Another verification gives $(\G\we\G)_N=\gen{e^*\we f^*,e^*\we g^*}$.\qed

(2) Groups $\At(\tl_r)$ and $C(\tl_r)$. For any ring $R$ with unity viewed as a right regular $R$-module and any right $R$-module $M$ the mapping $\lambda_M:M\to\text{Hom}_R(R,M)$ defined by $x.\lambda_M(m)=mx,x\in R$ is an $R$-isomorphism. Setting $M=R=R_3$ we have $\At(\tl_r)=\{\lambda_{R_3}(m)|m\in R_3\}$. Expand $m$ in the standard basis of $R_3$, $m=ue+qf+rg$. Then the matrix of $\phi=\lambda_{R_3}(m)$ is $\F=\begin{pmatrix}u&q&r\\0&u&q\\0&0&u\end{pmatrix}$. The matrices of mappings induced by $\phi^{-1}$ in $\G$ and $\G\we\G$ are $\F^{\Tr}$ and $\F^{\Tr}\we\F^{\Tr}$. Explicitly
\begin{equation}\label{inducedmatrices}\F^{\Tr}=\begin{pmatrix}u&0&0\\q&u&0\\r&q&u\end{pmatrix}\;\text{and}\;\F^{\Tr}\we\F^{\Tr}=\begin{pmatrix}u^2&0&0\\uq&u^2&0\\q^2-ur&uq&u^2\end{pmatrix}\end{equation}

We will show that $C(\tl_r)=\bu{\mathbb{Z}_p}$ by constructing a family of isomorphisms $\lambda_k:(G,\tl_r)\to(G,\tl_r^k)$ for every $k\in\bu{\mathbb{Z}_p}$. To this end, let us take $M=(R_3,\tl_r^k)$ and set $\lambda_k=\lambda_M(e)$. By definition of $\lambda_k$ we have
\begin{equation*} e.\lambda_k=e,\,f.\lambda_k=e(t^k-1),\,g.\lambda_k=e(t^k-1)^2\end{equation*}
Using the expansion $t^k-1=k(t-1)+\binom k2(t-1)^2\pmod {(t-1)^3}$ we conclude that $\Lambda_k=\begin{pmatrix}1&0&0\\0&k&\binom k2\\0&0&k^2\end{pmatrix}$ is the matrix of $\lambda_k$ in the standard basis. We shall need an explicit form of the associated matrices describing the action of $\lambda_k$ in $\G$ and $\G\we\G$, respectively. Put $l=k^{-1} \pmod p$ as usual. Then an easy calculation gives
\begin{equation}\label{asso1} (\Lambda_k^{-1})^{\Tr}=\begin{pmatrix}1&0&0\\0&l&0\\0&\binom l2&l^2\end{pmatrix},\end{equation}
\begin{equation}\label{asso2} (\Lambda_k^{-1})^{\Tr}\we(\Lambda_k^{-1})^{\Tr}=\begin{pmatrix}l&0&0\\\binom l2&l^2&0\\0&0&l^3\end{pmatrix}.\end{equation}

Unless stated otherwise we assume below that $p>3$. The degenerate case $p=3$ follows easily from the general one.

(3) Orbits of $\At(\tl_r)$ in $\x_p$. We identify $\x_p$ with $\mathbb{Z}_p^4$ via $x=ae^*+b_1e^*\we f^*+b_2e^*\we g^*+b_3g^*\we f^*\mapsto (a,b_1,b_2,b_3)$. We begin by listing all orbits in $\G^{C_p}$ and $\G\we\G$, respectively:  
\begin{align*} Z'_0&=\{(0)\},Z'_1=\{(a)|a\ne 0\},Z''_0=\{(0,0,0)\},\\Z''_{ij}&=\{(*,\ldots,*,\zeta^jb_i,0,\ldots,0)| b_i\in{\bu{\Z_p}}^2\},i=1,2,3;j=0,1\end{align*}
where the $*$ denotes an arbitrary element of $\Z_p$. For more complex orbits we need vectors $v_k(m)=(1,0,\ldots,m,0\ldots,0)\in\mathbb{Z}^4$ with the $m$ filling the $(k+1)$th slot, $k=1,2,3$ and running over $\bu{\Z_p}$.
\begin{Lem}\label{A'orbits} There are $3p+5$ orbits of $\At(\tl_r)$ in $\x_p$, namely
\begin{equation*}Z'_0\times Z''_0,\,Z'_1\times Z''_0,\,Z'_0\times Z''_{ij},\text{and}\;v_k(m)\At(\tl_r),k=1,2,3\end{equation*}\end{Lem}
\pf  The first two sets are clearly orbits. By \eqref{inducedmatrices} and every $i,j$

\noindent$(0,\ldots,\underset {i+1}{\zeta^j},0,\ldots,0).\phi=(0,*,\ldots,*,\zeta^ju^2,0\ldots,0)$ with the $*$ standing for an arbitrary element of $\mathbb{Z}_p$. This shows $Z'_0\times Z''_{ij}$ is the orbit of\\ $(0,\ldots,\underset{i+1}{\zeta^j},0,\ldots,0)$. Applying \eqref{inducedmatrices} again we have
\begin{equation}\label{v(m)phi} v_k(m).\phi=(u,*,\ldots,*,u^2m,0,\ldots,0)\end{equation}
From this one can see easily that $v_k(m)\At(\tl_r)$ has $(p-1)p^{k-1}$ elements. Another verification gives $v_k(m)\At(\tl_r)\cap v_k(n)\At(\tl_r)=\emptyset$ for $m\ne n$. Let us define $Z''_i=Z''_{i0}\cup Z''_{i1}$ and observe that $|Z''_i|=(p-1)p^{i-1}$ which gives $|Z'_1\times Z''_i|=(p-1)^2p^{i-1}$. Evidently $v_i(m)\in Z'_1\times Z''_i$ for all $m$ and therefore comparing cardinalities we arrive at the equality $Z'_1\times Z''_i=\bigcup_mv_i(m)\At(\tl_r)$. But clearly $\x_p=\bigcup Z'_l\times Z''_i,l=0,1;0\le i\le 3$ which completes the proof.\qed

(4) End of the proof.
\begin{Prop}\label{g'orbits} The nonzero orbits of $\g(\tl_r)$ in $\x_p$ are as follows:
\begin{equation*} Z'_0\times Z''_{ij},Z'_0\times Z''_2,\,Z'_1\times Z''_0,\,Z'_1\times Z''_2,\,Z'_1\times Z''_{3j},\,\text{and}\;v_1(m)\At(\tl_r),\end{equation*}
where $i=1,3,j=0,1$ and $m$ runs over ${\Z_p}^{\bullet}$.\end{Prop}
\pf By Corollary \ref{sizeoforbit} we need to determine the $\At(\tl_r)$-orbit containing $v\omega_k$ where $v$ runs over a set of generators of $\At(\tl_r)$-orbits of Lemma \ref{A'orbits}, and $\omega_k=\lambda_k\alpha_k^{-1},2\le k\le p-1$.

(i) For $\At(\tl_r)$-orbits $Z'_1\times Z''_0$ and $Z'_0\times Z''_{ij}$ generators are $e^*$ and $v_{ij}=(0,0,\ldots,\underset{i+1}\zeta^j,\ldots,0)$, respectively. In view of $e^*$ and $e^*\we f^*$ being fixed points for the action of $t$, and by \eqref{asso1},  \eqref{asso2} it is immediate that 
\begin{equation}\label{generators1}e^*\omega_k=le^*\,\text{and}\, v_{1j}\omega_k=l^2v_{1j},\end{equation}
hence $Z'_0\times Z''_{1j}$ and $Z'_1\times Z''_0$ are $\g(\tl_r)$-orbits. 

(ii) Next we take the generator $v_{20}=e^*\we g^*$. Noting that $(t-1)^2\bullet e^*\we g^*=0$, we use the expansion $\phi_l=l+\binom l2(t-1) \pmod{(t-1)^2}$ to derive
\begin{equation*}\phi_l\bullet e^*\we g^*=c e^*\we f^*+le^*\we g^*,\,c\in\Z_p.\end{equation*}
Applying $\lambda_k$ to the last equation we find with the help from \eqref{asso2}
\begin{equation}\label{generators2}e^*\we g^*\omega_k=c'e^*\we f^*+l^3e^*\we g^*,\,\text{for some}\,c'\in\Z_p.\end{equation}
The last equation shows that $v_{20}.\omega_k\in v_{21}\At(\tl_r)$ if $l$, hence $k$, is not a square, and $v_{20}.\omega_k\in v_{20}\At(\tl_r)$, otherwise. This means $v_{20}\g(\tl_r)=Z'_0\times (Z''_{20}\cup Z''_{21})=Z'_0\times Z''_2$ as needed.

The argument for the generator $v_{3j}=(0,0,0,\zeta^j)=\zeta^j f^*\we g^*$ of $Z_0'\times Z''_{3j}$ is almost identical.  Using the expansion $\phi_l=l+c_1(t-1)+c_2(t-1)^2 \pmod{(t-1)^3}$ we derive $\phi_l\bullet f^*\we g^*=(c_1+c_2)e^*\we f^*+c_1 e^*\we g^*+lf^*\we g^*$. Applying $\lambda_k$ we have by \eqref{asso2}
\begin{equation}\label{generators3}f^*\we g^*.\omega_k=c'_1e^*\we f^*+c_1l^2e^*\we g^*+l^4f^*\we g^*,\,c'_1,c_1\in\Z_p.\end{equation}
which shows $\zeta^j f^*\we g^*.\omega_k\in Z_{3j}$ for every $k$, hence $Z_0'\times Z''_{3j}$ is a $\g(\tl_r)$-orbit. 

(iii) We pause to mention that the above arguments settle the $p=3$-case. For, since $\x_3=\gen{e^*\we f^*,e^*\we g^*}$, by parts (i) and (ii) it has three nonzero orbits, namely $Z''_{1j},Z''_2,j=0,1$.

(iv) Here we take $v_1(m)=(1,m,0,0)$. Calculations in part (i) give $v_1(m).\omega_k=(l,l^2m,0,0)\in v_1(m)\At(\tl_r)$ by \eqref{v(m)phi}. That is, $v_1(m)\At(\tl_r)$ is a $\g(\tl_r)$-orbit for every $m\in\bu{\Z_p}$.

It remains to show that the last three sets of the Proposition are $\g(\tl_r)$-orbits.

(v) $Z'_1\times Z''_2$ is an orbit. By Lemma \ref{A'orbits} $Z'_1\times Z''_2=\bigcup_mv_2(m)\At(\tl_r)$ where $v_2(m)=e^*+me^*\we g^*$. Note that by \eqref{generators1} and \eqref{generators2} there holds $v_2(m).\omega_k=(l,c',l^3m,0)$. On the other hand we have by \eqref{inducedmatrices} $v_2(n).\phi=(u,uq,u^2n,0)$ where $u,q$ run over $\Z_p^{\bullet}$ and $\Z_p$, respectively. For every $l$ choosing $\phi=\phi(l,u^{-1}c',0)$ and $n=lm$ we obtain $v_2(m).\omega_k=v_2(n)\phi^{-1}$. Letting $k$ hence $l$ run over $\Z_p^{\bullet}$ we see that $v_2(m)\g(\tl_r)=\bigcup_nv_2(n)\At(\tl_r)$ which completes the proof.

(vi) Here we show that each $Z'_1\times Z''_{3j}$ is an orbit. By \eqref{generators1} and \eqref{generators3}
\begin{equation*}v_3(m).\omega_k=(1,0,0,m).\omega_k=(l,c',c'',ml^4)\,\text{for some}\,c',c''\in\Z_p.\end{equation*}
We seek  an $n$ such that
\begin{equation}\label{gvsA} v_3(m).\omega_k=v_3(n).\phi\,\text{for some}\,\phi\in\At(\tl_r).\end{equation}
By \eqref{inducedmatrices} $v_3(n).\phi=(u,q^2-ur,uq,u^2n)$ where $u,q,r$ take arbitrary values in $\Z_p^{\bullet}$ and $\Z_p$, respectively. Choosing $u,q,r$ such that

\noindent$u=l,q^2-ur=c',uq=c''$ and $n=ml^2$ fullfils \eqref{gvsA}. This yields the equality (*) $v_3(m)\g(\tl_r)=\displaystyle{\bigcup_{n\in m{\Z_p^{\bullet}}^2}}(1,0,0,n)\At(\tl_r)$. Therefore depending on $m\in{\bu{\mathbb{Z}_p}}^2$, or $m\notin{\bu{\mathbb{Z}_p}}^2$ the right hand side of (*) equals to $Z'_1\times Z''_{30}$ or $Z'_0\times Z''_{31}$, respectively. \qed

\subsubsection{$G=\Z_{p^e}\oplus\Z_p$}.\label{caseB} 

Our immediate goal is to classify nontrivial Hopf algebras in 

\noindent$\mathrm{Ext}(\k C_p,\k^{\Z_{p^2}\oplus \Z_p})$. We find it convenient to enlarge the scope of the problem by taking  $G=\Z_{p^e}\oplus\Z_p$ for {\em any} $e\ge 2$ as the amount of effort is the same as for $e=2$. As before our prime is odd, the even case is done in \cite{Ka1}. The end result is- 
\begin{Thm}\label{maincaseB} There are $2p+8$ distinct Hopf algebras in $\mathrm{Ext}(\k C_p,\k^G)$ if either $p>3$ or $e\ge 3$, and $16$ if $p=3$ and $e=2$.\end{Thm}
\pf We break up the proof in steps.
 
(1) Our first task is to describe the set of classes $[\tl]$ and their associated groups $\A,C(\tl)$. We need several preliminary observations. Every representation $\tl:C_p\to\Aut(G)$ is determined by $\tl(t)$. Let us write $\Gamma_e=\Aut(\Z_{p^e}\oplus\Z_p)$ and $\Gamma_e(p)$ for the set of all elements of order $p$ in $\Gamma_e$. It is clear that the mapping $\tl\mapsto\tl(t)$ sets up a bijection between the set $\{\tl\}$ and $\Gamma_e(p)$, and we will identify both sets. Furthermore the class $\text{eq}(\tl)$ of representations equivalent to $\tl$ corresponds to the ${\Gamma}_e$-conjugacy class of $\tl(t)$ denoted $\tl(t)^{\Gamma}$.  It follows that $[\tl]=\bigcup_{1\le k\le p-1}\tl(t^k)^{\Gamma}$.

$G$ has a natural basis $e_1,e_2$ comprised of generators of $\mathbb{Z}_{p^e},\mathbb{Z}_p$, respectively. Let $\epsilon$ be an endomorphism of $G$. We use the standard matrix representation of endomorphisms of direct sums to associate to $\epsilon$ a matrix $M(\epsilon)=\begin{pmatrix}a&\ov{b}\\cp^{e-1}&\ov{d}\end{pmatrix}$ relative to the basis $\{e_1,e_2\}$ with $a,b,c,d\in\mathbb{Z}_{p^e}$ and the bar over an $n\in\mathbb{Z}_{p^e}$ denoting the image of $n$ in $\mathbb{Z}_p$. The correspondence $\epsilon\mapsto M(\epsilon)$ extends to an isomorphism under the multiplication rule
\begin{equation*}\begin{pmatrix}a&\ov{b}\\cp^{e-1}&\ov{d}\end{pmatrix}\begin{pmatrix}a'&\ov{b'}\\c'p^{e-1}&\ov{d'}\end{pmatrix}=\begin{pmatrix}aa'+c'bp^{e-1}&\ov{ab'+bd'}\\(ca'+dc')p^{e-1}&\ov{dd'}\end{pmatrix}\end{equation*}
\begin{Lem}\label{Gamma}$\Gamma_e$ is the set of all matrices $\begin{pmatrix}a&\ov{b}\\cp^{e-1}&\ov{d}\end{pmatrix}$

\noindent satisfying $\ov{ad}\ne 0$\end{Lem}
\pf The natural epimorphism $G\to\mathbb{Z}_p\oplus\mathbb{Z}_p$ induces a homomorphism $\pi:\Gamma_e\to\Aut(\mathbb{Z}_p^2)$ via $\begin{pmatrix}a&\ov{b}\\cp^{e-1}&\ov{d}\end{pmatrix}\mapsto\begin{pmatrix}\ov{a}&\ov{b}\\0&\ov{d}\end{pmatrix}$. If $\gamma$ is invertible then so is $\pi(\gamma)$, and the latter is equivalent to $\ov{ad}\ne 0$. Conversely, if $\ov{ad}\ne 0$, then $a,d$ are units in $\mathbb{Z}_{p^e}$. One can check easily a factorization
\begin{equation}\label{factorization}\begin{pmatrix}a&\ov{b}\\cp^{e-1}&\ov{d}\end{pmatrix}=\begin{pmatrix}1&\ov{0}\\a^{-1}cp^{e-1}&\ov{1}\end{pmatrix}\begin{pmatrix}a&\ov{0}\\0&\ov{d}\end{pmatrix}\begin{pmatrix}a&\ov{a^{-1}b}\\0&\ov{d}\end{pmatrix}\end{equation}
which completes the proof.\qed

\begin{Lem}\label{Gamma(p)} {\rm(i)} $\Gamma_e(p)$ is the set of all matrices $\begin{pmatrix}1+ip^{e-1}&\ov{j}\\kp^{e-1}&\ov{1}\end{pmatrix}$;

{\rm(ii)} $|\Gamma_e(p)|=p^3$ regardless of $e$;

{\rm(iii)} $\Gamma_e(p)$ is a normal subgroup of $\Gamma_e$.
\end{Lem}
\pf (i) Assume $M=\begin{pmatrix}a&\ov{b}\\cp^{e-1}&\ov{d}\end{pmatrix}$ has order $p$. Then $\pi(M)$ has also order $p$ which implies $\ov{a}^p=\ov{1}=\ov{d}^p$, hence $\ov{d}=\ov{1}$ and $a\equiv 1\pmod p$. A straightforward induction on $r$ gives 
\begin{equation}\label{orderp} M^r=\begin{pmatrix}a^r+bc\binom{r}{2}p^{e-1}&\ov{rb}\\rcp^{e-1}&\ov{1}\end{pmatrix}\end{equation} 
whence $M^p=I$ iff $a^p=1$. But this condition on $a$ is equivalent to $a=1+ip^{e-1}$. (ii) and (iii) are easy consequences of (i). \qed

By the above Lemma $\Gamma_e(p)$ does not depend on $e$. We will omit $e$ from its notation below.  

\begin{Rem}\label{degenerationatp=2} All parts of this Lemma fail for $p=2$.\end{Rem}

\begin{Prop}\label{listofclasses} The set $\{[\tl]\}$ consists of five nontrivial elements.\end{Prop}
\pf (1) The first class of action is the one generated by $\tl_1$ with $\tl_1(t)=\text{diag}(1+p^{e-1},1)$, the diagonal matrix with entries $1+p^{e-1},\ov{1}$ on the main diagonal in that order. Using \eqref{factorization} one can see easily that the matrices $\text{diag}(1+kp^{e-1},\ov{1})$ form the center of $\Gamma(p)$. Since $\tl_1(t^k)=\text{diag}(1+kp^{e-1},\ov{1})$ it follows that $[\tl_1]=\{\tl_1^k|1\le k\le p-1\}$. As $\tl_1(t)$ is in the center $\At(\tl_1)=\Gamma_e,C(\tl_1)=\{1\}$ hence $\g(\tl_1)=\At(\tl_1)$. 

(2) Let $T_{\l}$ be the subset of lower triangular matrices in $\Gamma(p)$, $Z$ the center of $\Gamma(p)$ and $T'_{\l}=T_{\l}\setminus Z$. Fix one action $\tl_{\l}$ defined by $\tl_{\l}(t)=\begin{pmatrix}1&\ov{0}\\p^{e-1}&\ov{1}\end{pmatrix}$.

\begin{Lem}\label{lowertriangular} {\rm(i)} $T'_{\l}=\tl_{\l}^{\Gamma}$;

{\rm(ii)} $I(\tl_{\l},\tl_{\l}^k)\ne\emptyset$ for every $k$. In particular, $\text{diag}(1,\ov{k^{-1}})\in I(\tl_{\l},\tl_{\l}^k)$;

{\rm(iii)} $\At(\tl_{\l})=\left\{\begin{pmatrix}a&\ov{0}\\cp^{e-1}&\ov{a}\end{pmatrix}\right\}$ and $C(\tl_{\l})=A_p$.
\end{Lem}
\pf  (i) Pick another action $\tl$ with $\tl(t)=\begin{pmatrix}1+ip^{e-1}&\ov{0}\\jp^{e-1}&\ov{1}\end{pmatrix}, j\ne 0$. Matrices $\begin{pmatrix}1&\ov{0}\\kp^{e-1}&\ov{1}\end{pmatrix}$ lie in the center of $T_{\ell}$. By \eqref{factorization} $\tl_{\l}^{\Gamma}$ equals to $\{\tl_{\l}^{\gamma}\}$ where $\gamma$ runs over all upper triangular matrices in $\Gamma_e$. Choose a $\gamma=\begin{pmatrix}a&\ov{b}\\0&\ov{d}\end{pmatrix}$ and observe that $\gamma\in I(\tl_{\l},\tl)$ iff (*) $\tl_{\ell}(t)\gamma=\gamma\tl(t)$. One can see by a direct calculation that (*) holds iff 
\begin{align*}ai+bj&\equiv 0\pmod p\\a&\equiv jd\pmod p.\end{align*}
These congruences are equivalent to the conditions $b\equiv -aij^{-1}\!\!\pmod p$,

\noindent $d\equiv aj^{-1}\pmod p$ which gives (**) $I(\tl_{\l},\tl)=\{\begin{pmatrix}a&-\ov{aij^{-1}}\\cp^{e-1}&\ov{aj^{-1}}\end{pmatrix}\}$.

(ii) Take $\tl=\tl_{\l}^k$ and observe that $i=0,j=k$ for this action. Specifying $a=1, c=0$ in (**) yields (2).

(iii) Set $\tl=\tl_{\l}$ and note that $i=0,j=1$ in this case. Then (**) gives the assertion.\qed

(3) It remains to describe conjugacy classes in $\Gamma(p)\setminus T_{\l}$. Elements of this set are distinguished by the property-
\begin{Lem}\label{cyclicmodules} $\tl(t)\in\Gamma(p)\setminus T_{\l}$ iff the $C_p$-module $(G,\tl)$ is cyclic.\end{Lem}
\pf In one direction take $\tl(t)=\begin{pmatrix}1+ip^{e-1}&\ov{j}\\kp^{e-1}&\ov{1}\end{pmatrix}\in\Gamma(p)\setminus T_{\l}$. Then $\ov{j}\ne 0$ and therefore from $e_1\tl t=(1+ip^{e-1})e_1+je_2$ we have $e_2=j^{-1}e_1\tl (t-(1+ip^{e-1}))$ showing that $G$ is generated by $e_1$.

Conversely, assume $\ov{j}=0$. The subgroup $\gen{pe_1}$ is a $C_p$-submodule of $G$. Further, $G/\gen{pe_1}$ is a trivial $C_p$-module isomorphic to $\Z_p\oplus\Z_p$ which proves $(G,\tl)$ is not cyclic.\qed

We associate to an action $\tl\in\Gamma(p)\setminus T_{\l}$ with $\tl(t)=\begin{pmatrix}1+ip^{e-1}&\ov{j}\\kp^{e-1}&\ov{1}\end{pmatrix}$ the element $m(\tl)=jk$ of $\Z_{p^e}$. For an $n\in\Z_{p^e}$ we define $I(n)$ to be the ideal of $R$ generated by $p(t-1),(t-1)^2-np^{e-1}$ and $(t-1)^3$. $\ov{m}(\tl)$ is an invariant of $\tl(t)^{\Gamma}$ according to
\begin{Lem}\label{conjclassoftype2} {\rm(i)} In the foregoing notation $(G,\tl)\simeq R/I(m)$.

{\rm(ii)} Two actions $\tl,\tl'$ in $\Gamma(p)\setminus T_{\l}$ are equivalent iff $\ov{m}(\tl)=\ov{m}(\tl')$.\end{Lem}
\pf Let $R=\Z_{p^e}C_p$. Since $G=e_1R$ by the preceeding Lemma, both the assertions follow from the equality $I(m)=\text{ann}_Re_1$ for $m=m(\tl)$. In one direction, a simple calculation gives that $pe_1$ is a fixed point and $e_1\tl(t-1)^2=jkp^{e-1}e_1$. It follows that $e_1\tl g(t)=0$ for every generator $g(t)$ of $I(m)$ from the above list, whence $I(m)\subset\text{ann}_Re_1$. In the opposite direction we note every element of $R$ is congruent to some $n+m(t-1),n,m\in\Z_{p^e}$ modulo $I(m)$. Were $\text{ann}_Re_1\ne I(m)$, there would be an $n+m(t-1)$ with $e_1\tl(n+m(t-1))=0$, yet $n\ne 0$ or $m\not\equiv 0\pmod p$. But $e_1\tl(n+m(t-1))=(n+mip^{e-1})e_1+mje_2=0$ holds iff $m\equiv 0\pmod p$ and $n=0$ proving the equality in question.\qed

We single out three actions in $\Gamma(p)\setminus T_{\l}$,

\begin{equation}\label{threeactions}\tl^0=\begin{pmatrix}1&\ov{1}\\p^{e-1}&\ov{1}\end{pmatrix},\;\tl^1=\begin{pmatrix}1&\ov{1}\\p^{e-1}&\ov{1}\end{pmatrix},\;\tl^{\zeta}=\begin{pmatrix}1&\zeta\\p^{e-1}&\ov{1}\end{pmatrix}.\end{equation}
The next lemma completes the proof of the Proposition
\begin{Lem}\label{[]classes} $\Gamma(p)\setminus T_{\l}$ is the union of $[\tl^0]$, $[\tl^1]$ and $[\tl^{\zeta}]$.\end{Lem}
\pf By the formula \eqref{orderp} we have $m(\tl^r)=r^2m(\tl)$. The preceeding Lemma makes it clear that sets $[\tl^q],q=0,1,\zeta$ correspond to the orbits of ${\bu{\Z_p}}^2$ in $\Z_p$, namely $\{0\},{\bu{\Z_p}}^2,\zeta{\bu{\Z_p}}^2$.\qed

(4) We complete the proof of the main theorem of this section by computing the classifying groups and orbits for each of the five classes of actions. To begin with we select a basis for $\G$ dual to $\{e_i\}$ denoted by $\{e^*_i\}$. $C_p$ and $\Gamma_e$ act in $\G$ by \eqref{dualFaction} and $(f.\gamma)(g)=f(g\gamma^{-1}),f\in\G,g\in G$, respectively. These actions extend to $\mathrm{Alt}(G)=\G\we\G$ in the usual way. We note that $\mathrm{Alt}(G)$ is generated by $\beta=e_1^*\we e_2^*$ and the latter form has order $p$. For the future references we record
\begin{Lem}\label{contraposmatrices}{\rm(i)} Let $\begin{pmatrix}a&\ov{b}\\cp^{e-1}&\ov{d}\end{pmatrix}$ be the matrix of either $\gamma\in\Gamma$ or $t$ relative to $\{e_i\}$. The matrix of $\gamma^{-1}$ or $t$ relative to $\{e_i^*\}$ is $\begin{pmatrix}a&\ov{c}\\bp^{e-1}&\ov{d}\end{pmatrix}$
{\rm(ii)} There holds $\beta.\gamma^{-1}=ad\beta,\,t.\beta=\beta$, and $\mathrm{Alt}_N(G)=\mathrm{Alt}(G)$.\end{Lem}
\pf (i) is seen by a simple calculation. For (ii) we use part (i) to calculate $e_1^*\we e_2^*.\gamma^{-1}=(ae_1^*+ce_2^*)\we(bp^{e-1}e_1^*+de_2^*)=ade_1^*\we e_2^*$. Similartly $t.e_1^*\we e_2^*=ade_1^*\we e_2^*$. However in the case of $t$, $a=1+ip^{e-1}$ and $d=1$ by Lemma \ref{Gamma(p)}, which gives the second formula. Therefore $\phi_p(t).\beta=p\beta=0$ which proves the last assertion.\qed

(i) We take up the action $\tl_1$ of Proposition \ref{listofclasses}(1). 
\begin{Lem}\label{delta_1} $\mathrm{Ext}_{[\tl_1]}(\k C_p,\k^G)$ contains two distinct nontrivial Hopf algebras.\end{Lem}
\pf A simple calculation gives $\G^{C_p}=\gen{pe_1^*,e_2^*}$. As for $N(\G)$ we have $\phi_p(t).e_2^*=pe_2^*=0$ and $\phi_p(t).e_1^*=(\sum_{i=0}^{p-1}(1+p^{e-1})^i)e_1^*=pe_1^*$. It follows that $\G^{C_p}/N(\G)=\gen{\ov{e_2^*}}$ where $\ov{e_2^*}=e_2^*+N(\G)$. As noted in Proposition \ref{listofclasses}(1) $\g(\tl_1)=\At(\tl_1)=\Gamma_e$. By Lemma \ref{contraposmatrices} $\ov{e_2^*}.\gamma^{-1}=d\ov{e_2^*}$ and $\beta.\gamma^{-1}=ad\beta$. We conclude that $\x(\tl_1)\simeq \Z_p\oplus\Z_p$ with the action $(c_1,c_2).\gamma=(dc_1,adc_2)$. Now it is immediate that there are two nontrivial (i.e. $c_2\ne 0$) orbits, viz. $\{(0,c_2)\}$ and $\{c_1,c_2|c_1c_2\ne 0\}$.\qed

(ii) Next we consider $\tl_{\l}$ from Proposition \ref{listofclasses}(2).
\begin{Lem}\label{tll} There are $p+1$ distinct nontrivial Hopf algebras in $\mathrm{Ext}_{[\tl_{\l}]}(\k C_p,\k^G)$.\end{Lem}
\pf One can see easily with the help from Lemma \ref{contraposmatrices}  $\G^{C_p}=\gen{pe_1^*,e_2^*}$. Further $N(e_2^*)=pe_2^*=0$ and $N(e_1^*)=pe_1^*$. All in all we have $\G^{C_p}/N(\G)=\gen{\ov{e_2^*}}$ and $\x(\tl_{\l})=\gen{\ov{e_2^*},\beta}$. Using definition \eqref{derivedaction} we have $\ov{e_2^*}.\omega_k=(\phi_{k^{-1}}(t).\ov{e_2^*}).\lambda_k$ where $\lambda_k=\text{diag}(1,\k^{-1})$ by Lemma \ref{lowertriangular}. Since $\ov{e_2^*}$ is a fixed point, $\phi_{k^{-1}}(t).\ov{e_2^*}=k^{-1}\ov{e_2^*}$ and by Lemma \ref{contraposmatrices} $\ov{e_2^*}.\lambda_k=k\ov{e_2^*}$, hence $\ov{e_2^*}$ is fixed by $\omega_k$. A similar calculation gives $\beta.\omega_k=\beta$. Thus $\g(\tl_{\l})$-orbits coincide with $\At(\tl_{\l})$-orbits. For the latter we take $\phi\in\At(\tl_{\l})$ as in Lemma \ref{lowertriangular}(iii) and apply Lemma \ref{contraposmatrices} to get $\ov{e_2^*}.\phi^{-1}=\ov{a}\ov{e_2^*}$ and $\beta.\phi^{-1}=\ov{a}^2\beta$. It transpires that $\x(\tl_{\l})\simeq\Z_p^2$ with the action on the right by $(c_1,c_2).\phi^{-1}=(\ov{a}c_1,\ov{a}^2c_2)$. Now the argument in Proposition \ref{C_ptimesC_p} completes the proof.\qed 

(iii) Finally we tackle actions \eqref{threeactions}. We determine the groups 

\noindent$\At(\tl^q),C(\tl^q),q=0,1,\zeta$ and sets of intertwiners $\{\lambda_k|k\in C(\tl^q)\}$.
\begin{Lem}\label{intertwiners} {\rm(i)} $\At(\tl^q)=\left\{\begin{pmatrix}a&\ov{b}\\bqp^{e-1}&\ov{a}\end{pmatrix}\right\}$;

{\rm(ii)} $C(\tl^0)=A_p$ and for every $1\le k\le p-1$ $I(\tl^0,(\tl^0)^k)\ni\begin{pmatrix}1&\ov{0}\\0&k\end{pmatrix}$;

{\rm(iii)} If $q\ne 0$, then $C(\tl^q)\!=\!\{1,p-1\}$ and $I(\tl^q,(\tl^q)^{p-1})\!\ni\!\begin{pmatrix}1&\ov{0}\\qp^{e-1}&-\ov{1}\end{pmatrix}$
\end{Lem} 
\pf  (i) $\At(\tl^q)$ is the group of units of $\mathrm{End}_R(R/I(q))$. We pointed out in Theorem \ref{isoclassesAII}(2) that $\mathrm{End}_R(R/I(q))$ consists of mappings 

\noindent$\lambda(u):x\mapsto ux,u,x\in R/I(q)$. By Lemma \ref{conjclassoftype2}(i) $u=a\ov{1}+b\ov{(t-1)}$ where $\ov{r}=r+I(q)$ for $r\in R$. It is immediate that the matrix of $\lambda(u)$ relative to $\{\ov{1},\ov{(t-1)}\}$ is the one in part (i). 

(ii) and (iii) By Lemma \ref{conjclassoftype2} $C(\tl^q)=\{k|k^2q=q\}$. Clearly this formula implies $C(\tl^0)=A_p$ and $C(\tl^q)=\{1,p-1\}$ for $q\ne 0$. Let us write $\ov{R}=R/I(q)$ and denote by $\ov{R}^{(k)}$ the $C_p$-module $(\ov{R},(\tl^q)^k)$. By general principles for every $k\in C(\tl^q)$, $\text{Hom}_R(\ov{R},\ov{R}^{(k)})$ consists of mappings $\lambda(u),u\in\ov{R}$. Pick $\lambda(\ov{1})$ and observe that for every suitable $k$ the matrices of $\lambda(\ov{1})$ in the basis $\{\ov{1},\ov{t-1}\}$ are as given in (ii) and (iii), respectively.\qed

The last step of the proof of Theorems \ref{maincaseB} and \ref{mainp^4} is-
\begin{Lem}\label{classgroupsfortl^q} {\rm(i)} There are $p+1$ nontrivial distinct Hopf algebras in $\mathrm{Ext}_{[\tl^0]}(\k C_p,\k^G)$;

{\rm(ii)} There are two nontrivial distinct Hopf algebras in $\mathrm{Ext}_{[\tl^q]}(\k C_p,\k^G)$ for $q=1,\zeta$ if either $p>3$ or $e\ge 3$, and four otherwise.\end{Lem}
\pf (i) One can see easily that $\G^{C_p}(\tl^0)=\gen{e_1^*}$ and $N(\G(\tl^0))=pe_1^*$, hence $\G^{C_p}(\tl^0)/N(\G(\tl^0))=\gen{\ov{e_1^*}}$. By Lemma \ref{contraposmatrices}(ii) $\x(\tl^0)=\gen{\ov{e_1^*},\beta}$. Pick a $\gamma\in\At(\tl^0)$ as in Lemma \ref{intertwiners}. By Lemma \ref{contraposmatrices} there holds
$\ov{e_1^*}.\gamma^{-1}=\ov{a}\ov{e_1^*}\,\text{and}\,\beta.\gamma^{-1}=\ov{a}^2\beta$. This type of action occured in Proposition \ref{C_ptimesC_p} whose argument yields $p+1$ nontrivial $\At(\tl^0)$-orbits. Turning to $\g(\tl^0)$-orbits, pick a $\lambda_k=\text{diag}(1,k)$ from the preceeding lemma. Since $\ov{e_1^*},\beta$ are fixed by $t$ we have $\ov{e_1^*}.\omega_k=(\phi_{k^{-1}}.\ov{e_1^*}).\lambda_k=k^{-1}\ov{e_1^*}$ and $\beta.\omega_k=(\phi_{k^{-1}}.\beta).\lambda_k=k^{-2}\beta$. This shows that $\g(\tl^0)$-orbits coincide with $\At(\tl^0)$ ones, and the proof is complete.

(ii) A straighforward calculation gives $\G^{C_p}(\tl^q)=\gen{pe_1^*}$. For calculation of $N(\G(\tl^q))$ we employ \eqref{orderp} which gives readily that
\begin{equation*}\phi_p(t).e_1^*=[\sum_{r=0}^{p-1}(1+q\binom{r}{2}p^{e-1})]e_1^*+(\sum_{r=0}^{p-1}r)e_2^*\end{equation*}
As $\sum_{r=0}^{p-1}r=\binom{p}{2}$ and $pe_2^*=0$ we conclude 

\noindent$\phi_p(t).e_1^*=(p+q(\sum_{r=0}^{p-1}\binom{r}{2})p^{e-1})e_1^*$. Similarly one can derive
\begin{equation*}\phi_p(t).e_2^*=q(\sum_{r=0}^{p-1}r)p^{e-1}e_1^*+pe_2^*=0\end{equation*}
Next we note that an elementary calculation gives $\sum_{r=0}^{p-1}\binom{r}{2}=\binom{p}{3}$. Let us put $c(p)=p+q\binom{p}{3}p^{e-1}$. We observe that if $p>3$, then $c(p)\equiv p\pmod {p^e}$. For $p=3$ and either $e\ge 3$ or $e=2$ and $q=1$, $c(3)=3u$, where $u$ is a unit in $\Z_{p^e}$. In the exceptional case $e=2$ and $q=2$, $c(3)=9$. This translates into $\phi_p(t).e_1^*=pue_1^*$ for all $p,e,q$, except for the exceptional case where $\phi_3(t).e_1^*=0$. We conclude that $N(\G(\tl^q))=\gen{pe_1^*}$ in the regular case and it is zero, otherwise. In consequence 
\begin{align*}\x(\tl^1)&=\gen{\beta}\;\text{for all}\,p,e\\\x(\tl^{\zeta})&=\gen{\beta}\;\text{if}\,p>3\,\text{or}\,e\ge 3\\\x(\tl^2)&=\gen{3e_1^*,\beta}\;\text{if}\,p=2=e.\end{align*}
By Lemmas \ref{contraposmatrices}(ii),\ \ref{intertwiners}(i) $\beta.\phi^{-1}=\ov{a}^2\beta$ for every $\phi\in\At(\tl^q)$. It follows that there are two nontrivial $\At(\tl^q)$-orbits in $\x(\tl^q)$ in the regular case and also for $\x(\tl^1)$ in all cases, namely $\{cq\beta|c\in{\bu{\Z_p}}^2$ for $q=1,\zeta\}$. Using Lemma \ref{intertwiners}(iii) it is immediate that $\beta.\omega_{p-1}=\beta$. That says $\g(\tl^q)$-orbits coinside with $\At(\tl^q)$-orbits. In the exceptional case $3e_1^*.\phi^{-1}=\ov{a}(3e_1^*)$ and $3e_1^*.\omega_2=-3e_1^*$. It follows that $\At(\tl^2)$ and $\g(\tl^2)$ act on $\x(\tl^2)$ by $(c_1,c_2).\phi^{-1}=(\ov{a}c_1,\ov{a}^2c_2)$ and  $(c_1,c_2).\omega_2=(-c_1,c_2)$ with the usual identification $\x(\tl^2)\simeq{\Z_3}^2$. By the argument of Proposition \ref{C_ptimesC_p} there are four nontrivial $\At(\tl^2)$-orbits. One can check directly that the mapping $(c_1,c_2)\mapsto(-c_1,c_2)$ preserves the orbits, completing the proof.\qed


\section{Some old classification results revisited}

The first result concerns the G. Kac's $8$-dimensional Hopf algebra \cite{{KP},{Mas2}} which we denote by $H_8$.
\begin{Thm}\label{KPM} There is a unique semisimple, nontrivial $8$-dimensional Hopf algebra.\end{Thm}
\pf It is easy to see that every Hopf algebra $H$ as in the Theorem is isomorphic to $\k^4\oplus M_2(\k)$ as algebra where $M_2(\k)$ is the algebra of $2\times 2$ matrices. Applying this remark to $H^*$ we conclude that $H^*$ has exactly $4$ characters, hence $G(H)$ has order $4$. Thus $H$ is almost abelian, hence $H\in\mathrm{Ext}(\k C_2,\k^{G(H)})$. By Theorem \ref{mainthm}(II) the number of nontrivial isotypes in $\mathrm{Ext}_{[\tl]}(\k C_2,\k^G)$ equals to the number of nontrivial $\A$-orbits in $H^2_c(\tl)$ for every action $\tl$ of $C_2$ on $G$. By Corollary \ref{basicisic} that number coincides with the number of nontrivial $\A$-orbits in $\X$. For every cyclic group $C_n$, $\mathrm{Alt}(C_n)$ is trivial. Hence, were $G=C_4$ we would have $\X=\G^{C_2}/N(\G)$ by Lemma \ref{gisomorphism}(ii) and therefore $\X$ does not have nontrivial orbits. We take up the remaining case $G=G(H)=C_2\times C_2$. Let $\{x_1,x_2\}$ be a basis for $G$ and $\{x_1^*,x_2^*\}$ its dual. There is only one equivalence class of actions on $G$. We choose the action $x_1\tl t=x_2,x_2\tl t=x_1$. A routine verification gives $\G^{C_2}=N(\G)=\gen{x_1^*x_2^*}$. Thus by Lemma \ref{gisomorphism} $\X\simeq\un{a}(Z^2_N(\tl))$ and by Proposition \ref{hopf2cohomology}(3) we have $\un{a}(Z^2_N(\tl))=\mathrm{Alt}_N(G)$. Further, it is immediate that $\mathrm{Alt}_N(G)=\mathrm{Alt}(G)$ and the latter consists of one nonzero element. This shows that $\X$ has one nontrivial $\A$-orbit, and the proof is complete.\qed

With a small additional effort one can give a presentation of $H_8$ by generators and relations. For two vectors $a=x_1^{j_1}x_2^{j_2},b=x_1^{k_1}x_2^{k_2}$ we let $\det(a,b)=j_1k_2-j_2k_1$.
\begin{Prop}\label{defrelnsH8} $H_8$ is generated as algebra by $x_1^*,x_2^*,t$ subject to the relations
\begin{align*}{x_1^*}^2&={x_2^*}^2=t^2=1\\  
tx_1^*t^{-1}&=x_2^*,tx_2^*t^{-1}=x_1^*\end{align*}
The coalgebra structure is specified by
\begin{equation*}\Delta(t)=(\sum_{a,b\in G}\iota^{-\det(a,b)}p_a\o p_b)t\o t,\,\text{where}\,\iota^2=-1.\end{equation*}
In addition the equations $S(x_i^*)=x_i^*, i=1,2, S(t)=t$ and $\epsilon(x_1^*)=\epsilon(x_2^*)=\epsilon(t)=1$ determine the antipode and augmentation.
\end{Prop}
\pf Since $H_8$ is a special cocentral extensions $H_8=\k\G\#\k C_2$ as algebra. With $t$ a generator of $C_2$ the algebra relations follow immediately. By \eqref{DeltaH} 

\noindent$\Delta(t)=\displaystyle{(\sum_{a,b\in G}\tau(t,a,b)p_a\o p_b)t\o t}$ where $\tau(t,a,b)\in\X$. As $\X$ has only one nonzero element, the latter provided by Proposition \ref{hopf2cohomology}(3ii), we have $\tau(t,a,b)=s_{1,2}\delta g$. A straightforward calculation gives 

\noindent$\tau(t,a,b)=\iota^{-\det(a,b)}$.

We find the antipode by using \cite[Prop. 4.7]{M}. In our case, i.e. for a special cocentral extension, the formula specializes to $S(p_at)=\tau^{-1}(t,a^{-1},a)p_{a^{-1}\tl t}t^{-1}$. Since $a^2=1$ and $\tau(t,a,a)=1$, we obtain $S(t)=\sum_{a}S(p_at)=\sum_{a}p_{a\tl t}t=t$. The rest of the Proposition is self-evident.\qed

A. Masuoka \cite{Mas2} presents $H_8$ by a different set of generators and relations. The two are related by replacing $t$ with $z=gx_1^*t$. The set $\{x_1^*,x_2^*,z\}$ generates $H_8$ and one can derive all relations of \cite[Thm. 2.13]{Mas2}, with one exception, viz. $S(z)=\frac{1}{2}(-\epsilon+x_1^*+x_2^*+x_1^*x_2^*)z$. We leave the details to the reader.

We take up the problem of classifying isotypes of Hopf algebras $H$ of dimension $2n^2$ with $G(H)=\mathbb{Z}_n\times\mathbb{Z}_n$ for an odd $n$. Put differently we want to determine the isotypes of $\mathrm{Ext}(\k C_2,\k^{\mathbb{Z}_n\times\mathbb{Z}_n})$. We let $G=\mathbb{Z}_n\times\mathbb{Z}_n$

Following the general procedure we split up the argument into steps. 

(1) A survey of actions.

We will assume $n=p_1^{e_1}\cdots p_m^{e_m}$ is the prime decomposition of $n$. We let $G(i)$ denote the $p_i$-primary summand of $G$. Clearly $G(i)=\mathbb{Z}_{p_i^{e_i}}\oplus\mathbb{Z}_{p_i^{e_i}}$ and $G=\oplus G(i)$. Every $G(i)$ is invariant under any automorphism of $G$, in particular under any action of $C_2$. Since every $p_i$ is odd
$\mathbb{Z}_{p_i^{e_i}}C_2=\mathbb{Z}_{p_i^{e_i}}\epsilon_0\oplus\mathbb{Z}_{p_i^{e_i}}\epsilon_{-1}$ where $\epsilon_0=\frac{1+t}{2},\,\epsilon_{-1}=\frac{1-t}{2}$. Idempotents $\epsilon_{\nu}$ induce a splitting $G(i)=G(i)\epsilon_0\oplus G(i)\epsilon_{-1}$ into a direct sum of subgroups on which $t$ acts as $\pm\text{id}$. Therefore for every action $\tl$ we can write $G$ as
\begin{equation}\label{Gdecomposition}G=G_0\oplus G_{-1}\oplus G_{0,-1},\,\text{where}\end{equation} 
\begin{align*} &G_0=\oplus\{G(i)|\;t|_{G(i)}=\text{id}\},\hspace{.2in} G_{-1}=\oplus\{G(i)|\;t|_{G(i)}=-\text{id}\},\,\text{and}\\&G_{0,-1}=\oplus\{G(i)|\; t|_{G(i)}\ne\pm\text{id}\}.\end{align*}
Every equivalence class of actions is determined by its decomposition \eqref{Gdecomposition}. 

(2) Classifying groups.

First we show that $\G^{C_2}/N(\G)=(0)$. Pick $\chi\G^{C_2}$. Then $N(\chi):=(1+t).\chi=2\chi$. Since $2$ is a unit in $\mathbb{Z}_n$, $\chi\in N(\g)$, which proves our assertion. By Lemma \ref{gisomorphism}(iii) $\X=\mathrm{Alt}_N(G)$. Consider an alternate mapping $\beta:G\times G\to\mathbb{Z}_n$. It is apparent that $\beta(g,h)=0$ whenever $g,h$ lie in different components $(G(i)$ of decomposition \eqref{Gdecomposition}. For $g,h\in G_0$ $(1+t).\beta(g,h)=2\beta(g,h)$ and similarly for if $g,h\in G_{-1}$. It transpires that $(1+t).\beta(g,h)=0$ iff $\beta(g,h)=0$ for every $\beta:G_{\nu}\times G_{\nu}\to\mathbb{Z}_n,\,\nu=0,-1$. We conclude that $\X=0$ if $G_{0,-1}=0$. 

The above discussion shows that $\mathrm{Alt}_N(G)=\mathrm{Alt}_N(G_{0,-1})$. Let us renumber the prime divisors of $n$ so that $G_{0,-1}=\oplus_{i=1}^r\,G(i)$. We noted above that $G(i)=G(i)\epsilon_0\oplus G(i)\epsilon_{-1}$ and since $\mathbb{Z}_{p_i^{e_i}}$ is an indecomposable group, $G(i)\epsilon_{\nu}\simeq\mathbb{Z}_{p_i^{e_i}}$. Therefore we can select a basis $\{a_i,b_i\}$ of $G(i)$ with $a_i,b_i$ generating $G(i)\epsilon_0,G(i)\epsilon_{-1}$, respectively and both of order $p_i^{e_i}$. Set $a=\sum a_i,\,b=\sum b_i$ and observe that $a,b$ generate subgroups $G_{0,-1}\epsilon_{\nu}, \nu=0,-1$, respectively. Let us write $n(\tl)=\prod\{p_i^{e_i}|t|_{G(i)}\not=\pm\text{id}\}$. Set $a=\sum a_i,\,b=\sum b_i$ and observe that $a,b$ generate subgroups $G_{0,-1}\epsilon_{\nu}, \nu=0,-1$, respectively. In addition both subgroups $\gen{a},\gen{b}$ are cyclic of order $n(\tl)$, hence $G_{0,-1}\simeq\mathbb{Z}_{n(\tl)}\times\mathbb{Z}_{n(\tl)}$. It follows that $\mathrm{Alt}(G_{0,-1})$ is cyclic on a generator, say, $\beta_0$ defined by $\beta_0(a,b)=1_{\mathbb{Z}_{n(\tl)}}$. The calculation $(1+t).\beta_0(a,b)=\beta_0(a,b)+\beta(a,-b)=0$ gives the equality $\mathrm{Alt}_N(G_{0,-1})=\mathrm{Alt}(G_{0,-1})$. It follows that $\X=\mathrm{Alt}(G_{0,-1})\simeq \mathbb{Z}_{n(\tl)}$. We observe that since $\X\simeq H^2_c(\k C_2,\k^G,\tl)$ that formula implies a result of A. Masuoka \cite[Thm. 2.1]{Mas5} on $\mathrm{Opext}(\k C_2,\k^G)$.

We summarize
\begin{Thm}\label{isotypesforCn^2} {\rm(1)} If $\tl$ is such that $G_{0,-1}=0$, then $\mathrm{Ext}_{[\tl]}(\k C_2,\k^G)$ has a unique Hopf algebra $\k[G\rtimes C_2]$ where $G\rtimes C_2$ is the semidirect product with respect to $\tl$.

{\rm(2)} For $\tl$ with a nonzero $G_{0,-1}$ the isotypes in $\mathrm{Ext}_{[\tl]}(\k C_2,\k^G)$ correspond bijectively to the subgroups of $\mathbb{Z}_{n(\tl)}$. The trivial subgroup of $\mathbb{Z}_{n(\tl)}$ corresponds to a unique trivial Hopf algebra $\k[G\rtimes C_2]$.\end{Thm}
\pf It remains to compute the orbits of $\A$ in $\mathrm{Alt}(G_{0,-1})$. First off, every $\phi\in\A$ preserves $G_{0,-1}\epsilon_{\nu}$, whence 
$$a\phi=ua,\;b\phi=vb\;\text{for some}\;u,v\in\mathbb{Z}_{n(\tl)}^{\bullet}.$$
Therefore $(\beta.\phi)(a,b):=\beta(a\phi^{-1},b.\phi^{-1})=u^{-1}v^{-1}\beta(a,b)$. This shows that transwering action of $\A$ along the isomorphism $\beta\mapsto\beta(a,b):\mathrm{Alt}(G_{0,-1})\overset\sim\to\mathbb{Z}_{n(\tl)}$ we get the action $m.\phi=u^{-1}v^{-1}m,\,m\in\mathbb{Z}_{n(\tl)}$. It becomes clear that orbits are exactly sets of generators of cyclic subgroups of $\mathbb{Z}_{n(\tl)}$, which completes the proof.\qed
 
\vspace{.1in}
\section{Appendices}

{\bf Appendix 1: Crossed product splitting of abelian extensions}

\begin{Prop}. Let $H$ be an extension of $\k F$ by $\k^G$. Then $H$ is a crossed product of $\k F$ over $\k^G$.\end{Prop}

\pf First observe that $H$ is a Hopf-Galois extension of $\k^G$ by $\k F$ via $\rho_{\pi}=(\text{id}\otimes\pi)\Delta_H:H\to H\otimes\k F$, see e.g. the proof of \cite[3.4.3]{Mo}, hence by \cite[8.1.7]{Mo} $H$ is a strongly $F$-graded algebra. Setting $H_x=\{h\in H|\rho_{\pi}(h)=h\otimes x\}$ we have $H=\displaystyle{\oplus_{x\in F}}H_x$ with $H_1=\k^G$ and $H_xH_{x^{-1}}=\k^G$ for all $x\in F$. Next for every $a\in G$ we construct elements $u(a)\in H_x,\,v(a)\in H_{x^{-1}}$ such that 
\begin{align*} &u(a)v(a)=p_a,\,p_au(a)=u(a),v(a)p_a=v(a),\,\text{and}\\&u(a)v(b)=0\;\text{for all}\,a\ne b.\end{align*}
Indeed, were all $uv,u\in H_x,v\in H_{x^{-1}}$ lie in $\text{span}\{p_b|b\ne a\}$, then so would $H_xH_{x^{-1}}$, a contradiction. Therefore for every $a\in G$ there are $u\in H_x,v\in H_{x^{-1}}$ such that $uv=\sum c_bp_b,c_a\ne 0$. Setting $u(a)=\displaystyle\frac{1}{c_a}p_au,v(a)=vp_a$ we get elements satisfying the first three properties stated above. Furthermore, the last property also holds because $u(a)v(b)=p_au(a)v(b)p_b=p_ap_bu(a)v(b)=0$. It follows that the elements $u_x=\sum_{a\in G}u(a),v_x=\sum_{a\in G}v(a)$ satisfy $u_xv_x=1$ hence, as $H$ is finite-dimensional, $v_xu_x=1$ as well. Thus $u_x$ is a $2$-sided unit in $H_x$. 

Now define $\gamma:\k F\to H$ by $\gamma(x)=\displaystyle\frac{1}{\epsilon_H(u_x)}u_x$. One can see immediately that $\gamma$ is a convolution invertible mapping satisfying $\rho_{\pi}\circ\gamma=\gamma\o\text{id},\gamma(1_F)=1$ and $\epsilon_H\circ\gamma=\epsilon_F$. Thus $\gamma$ is a section of $\k F$ in $H$, which completes the proof.\qed
\vspace{.1in}

{\bf Appendix 2: Non-splitting of $\X$ as $\A$-module for $p=2$}
\vspace{.1in}

We take a closer look at the exact sequence $\G^{C_p}/N(\G)\rightarrowtail\X\twoheadrightarrow\un{a}(Z^2_N(\tl))$ of Lemma \ref{gisomorphism}. We know by Theorem \ref{hopf2cohomology} that for $p>2$ $\un{a}(Z^2_N(\tl))=\mathrm{Alt}_N(G)$ and the above sequence splits up, that is $\X\simeq\G^{C_p}/N(\G)\times\mathrm{Alt}_N(G)$ as $\A$-modules. We want to show that this is not the case for $p=2$.

Let $G$ be an elementary $2$-group of rank $n$ and $\tl$ be the trivial action. By the argument of part (2) of Proposition \ref{hopf2cohomology} our assumptions imply $\G^{C_p}/N(\G)=\G$ and $\un{a}(Z^2_N(\tl))=\mathrm{Alt}(G)$. The main result of this Appendix is
\begin{Thm}\label{nonsplittingX} Let $G$ be a $2$-elementary group of rank $n>2$. The sequence of $\At(\t)$-modules 
\begin{equation*}\G\to\x(\t)\to\mathrm{Alt}(G)\end{equation*}
does not split.\end{Thm}
\pf Will be given in steps. To simplify notation we write $\x$ and $\At$ for $\x(\t)$ and $\At(\t)$.

(1) Let $S$ be a copy of $\mathrm{Alt}(G)$ in $Z^2(G,\k^{\b})$ constructed in Proposition \ref{hopf2cohomology}(2). Clearly $S\subset Z^2_N(\t)$ and complements $B^2_N(\t)$. Passing on to $\x$ the image of $S$, denoted by $S$, forms a complement to $\G$. Fix a basis $\{x_i|1\le i\le n\}$ of $G$ and let $\{x_i^*|1\le i\le n\}$ be its dual in $\G$. Observe that $\Phi:B^2_N(\t)\to\G$ acts in the present case by $\Phi(\delta f)=f^2$. Let $b_i:G\times G\to\k^{\b}$ be the bimultiplicative map defined by 
\begin{equation*}b_i(x_i,x_i)=-1,b_i(x_k,x_l)=1\;\text{for }\;(k,l)\ne(i,i).\end{equation*}
\begin{Lem}\label{Phi}\rm{(1)} $\Phi(b_i)=x_i^*$ for all $i$;\\
\rm{(2)} $\mathrm{Alt}(G)\subset\ker\Phi$.
\end{Lem}
\pf (1) Recall $B^2(G,\k^{\b})$ is the subgroup of all symmetric functions of $Z^2(G,\k^{\b})$, hence $b_i\in B^2(G,\k^{\b})$ and therefore $b_i=\delta f_i$ for some $f_i:G\to\k^{\b}$. Then
\begin{equation*} b_i(x_j,x_j)=\delta f_i(x_j,x_j)=f_i(x_j)f_i(x_j)f_i(x_j^2)^{-1}=f_i^2(x_j).\end{equation*}
We note that as $b_i^2=\epsilon$, $b_i$ lies in $B^2_N(G,\k^{\b})$, hence $f_i^2\in\G$ and as $f_i^2(x_j)=(-1)^{\delta_{ij}}$ $f_i^2=x_i^*$. This proves (1). 

(2) Elements of $\mathrm{Alt}(G)$ are symmetric functions, hence 

\noindent$\mathrm{Alt}(G)\subset B^2(G,\k^{\b})$. By part (1) for every $\alpha=\delta f\in\mathrm{Alt}(G)$ $\Phi(\alpha)=f^2=\epsilon$ as $\alpha(x,x)=1$.\qed
 
(2) Let $\G\we\G$ be the exterior square of $\G$. There is a well-known identification $\mathrm{Alt}(G)=\G\we\G$. In the additive notation $\G\we\G$ has a standard basis $x_i^*\we x_j^*$ where $x_i^*\we x_j^*(x_k,x_l)=\delta_{ik}\delta_{jl}$. Passing on to $S$ we write $s_{x_i^*\we x_j^*}$ as $s_{\gen{i,j}}$ which by the definition of $s_{\alpha}$ is given by
$$s_{\gen{i,j}}(x_k,x_l)=\begin{cases}1,&\text{if}\;\{k,l\}=\{i,j\}\;\text{and}\;k<l\\0,&\text{else}.\end{cases}$$
We note the equality $s_{\gen{i,j}}=s_{\gen{j,i}}$. 
Pick $\phi\in\At$ and let $\phi^*:\G\to\G$ be the transpose of $\phi$, i.e. $(\chi.\phi^*)(g)=\chi(g.\phi),\chi\in\G,g\in G$. If $M(\phi)$ is the matrix of $\phi$ in the basis $\{x_k\}$ then $M(\phi^*)=M(\phi)^{\text{tr}}$ is the matrix of $\phi^*$ in the dual basis. Therefore the matrix of the mapping $\widehat{\phi}$, $(\chi.\widehat{\phi})(g)=\chi(g.\phi^{-1})$ induced by $\phi$ in $\G$ is $M(\phi^{-1})^{\text{tr}}$. Next we describe action of $\At$ in $\x$
\begin{Lem}\label{actionofA} Suppose $\phi\in\At$ and $M(\phi^{-1})=(a_{kl})$. $\phi$ acts in $\x$ as follows
\begin{equation}\label{Aaction}s_{\gen{i,j}}.\phi=s_{x_i^*\we x_j^*.\phi}+\sum\limits_{k=1}^na_{ki}a_{kj}x_k^*.\end{equation}
\end{Lem}
\pf One can see easily that the mapping $\un{a}$ is $\At$-linear therefore $\un{a}(s_{\gen{i,j}}.\phi)={x_i^*\we x_j^*.\phi}$. We also know $\un{a}(s_{\alpha})=\alpha$ for every $\alpha\in\mathrm{Alt}(G)$ which gives
\begin{equation}\label{intermediary}s_{\gen{i,j}}.\phi=s_{x_i^*\we x_j^*.\phi}+\ov{b},.\end{equation}
where $\ov{b}:=b\ker\F\in B^2_N(\t)/\ker\F$. By Lemma \ref{Phi} the set $\{\ov{b}_k\}$ forms a basis for $B^2_N(\t)/{\ker\Phi}$, hence $\ov{b}=\sum\limits_{k=1}^nc_k\ov{b}_k,\,c_k\in\Z_2$. Since $s_{\alpha}(x_k,x_k)=0$ for every $\alpha$ and $k$, evaluating \eqref{intermediary} at $(x_k,x_k)$ yields
\begin{align*}c_k&=s_{\gen{i,j}}.\phi(x_k,x_k)=s_{\gen{i,j}}(x_k\phi^{-1},x_k\phi^{-1})\\&=s_{\gen{i,j}}(\sum\limits_ia_{ki}x_i,\sum\limits_ja_{kj}x_j)=a_{ki}a_{kj}\end{align*}
A reference to Lemma \ref{Phi}(1) completes the proof.\qed

It is well known that $\At$ is generated by transvections, linear mappings $t_{pq}:x_p\to x_p+x_q,x_r\to x_r,r\ne p$. Since $t_{pq}^{-1}=t_{pq}$ and the matrix of $t_{pq}^*$ is $M(t_{pq})^{\text{tr}}$ we have readily 
\begin{align*}x_k^*.t_{pq}&=x_k,k\ne q,\\x_q^*.t_{pq}&=x_q^*+x_p^*.\end{align*}
We see that $t_{pq}$ induces the transvection $t_{qp}$ in $\G$. In consequence we have
\begin{Lem}\label{transvectionsinAlt} Action of transvections on the standard basis of $\mathrm{Alt}(G)$ is given by
\begin{align*}x_i^*\we x_j^*.t_{pq}&=x_i^*\we x_j^*\;\text{if}\;q\ne i,j\;\text{or}\;(p,q)=(i,j),(j,i)\\
              x_i^*\we x_j^*.t_{pi}&=x_i^*\we x_j^* +x_p^*\we x_j^*,\,p\ne j\\
              x_i^*\we x_j^*.t_{pj}&=x_i^*\we x_j^* +x_i^*\we x_p^*,\,p\ne i.\qed\end{align*}\end{Lem}
With the help of Lemma \ref{actionofA} we deduce
\begin{Lem}\label{transvectionsinS} Action of transvections on generators of $S$ is given by
\begin{align*} s_{\gen{i,j}}.t_{pq}&= s_{\gen{i,j}},\;\text{if}\;q\ne i,j\\
               s_{\gen{i,j}}.t_{ij}&=s_{\gen{i,j}}+x^*_i\\
               s_{\gen{i,j}}.t_{ji}&=s_{\gen{i,j}}+x^*_j\\
               s_{\gen{i,j}}.t_{pi}&=s_{\gen{i,j}}+s_{\gen{p,j}}\\
               s_{\gen{i,j}}.t_{pj}&=s_{\gen{i,j}}+s_{\gen{i,p}}.\end{align*}
               \end{Lem}             
\pf In view of Lemmas \ref{actionofA} and \ref{transvectionsinAlt} we need only to calculate the $\G$- components. If $(p,q)\ne (i,j),(j,i)$, then for the entries  of $M(t_{pq})$ there holds $a_{ki}=0$ or $a_{kj}=0$ for every $k$. In $M(t_{ij}),M(t_{ji})$ we have $a_{ki}a_{kj}=1$ only for $k=i,j$, respectively.\qed

(3) End of the Proof. Suppose there is an $\At$-linear section $\zeta:\mathrm{Alt}(G)\to\x$ splitting $\underline{a}$. Say 
\begin{equation}\label{Asection}\zeta(x_i^*\we x_j^*)=\chi_{ij}+s_{\gen{i,j}},\,\chi_{ij}\in\G.\end{equation}      
Then there holds
\begin{equation}\label{Alinear} \zeta(x_i^*\we x_j^*.t_{pq})=(\chi_{ij}+s_{\gen{i,j}}).t_{pq}\;\text{for all}\; p,q.\end{equation}
Let us expand $\chi_{ij}$ in the basis $\{x_k^*\}$, 
\begin{equation*} \chi_{ij}=\sum\limits_kc^{ij}_kx_k^*.\end{equation*}  
Observe the equality $\chi_{ij}.t_{pq}=\chi_{ij}+c^{ij}_qx_p^*$. Next specialize \eqref{Alinear} to $p=i,q=j$ or $p=j,q=i$. Then Lemmas \ref{transvectionsinAlt} and \ref{transvectionsinS} give $c_j^{ij}x_i^*+x_i^*=0$ and $c_i^{ij}x^*_j+x_j^*=0$, respectively. We see that $c_i^{ij}=c_j^{ij}=1$, that is $\chi_{ij}=x_i^*+x_j^*+\sum\limits_{k\ne i,j}c_k^{ij}x_k^*$.  Note that if $n=2$ we have shown that $\Z_2(x_1^*+x_2^*+s_{\gen{1,2}})$ is an $\At$-complement to $\G$. Suppose $n>2$. For every $q\ne i,j$ we have by \eqref{Alinear} and Lemmas \ref{transvectionsinAlt} and \ref{transvectionsinS} the equality
\begin{equation*} \chi_{ij}+s_{\gen{i,j}}=\chi_{ij}.t_{iq}+ s_{\gen{i,j}}.t_{iq}\end{equation*}
Using $\chi_{ij}.t_{iq}=\chi_{ij}+c_q^{ij}x_i^*$ and $s_{\gen{i,j}}.t_{iq}=s_{\gen{i,j}}$ we conclude $c_q^{ij}=0$. Thus $\chi_{ij}=x^*_i+x^*_j$ for all $i,j$. 

Next pick $p\ne i,j$, and apply \eqref{Alinear}. We have 
\begin{equation*}\zeta(x_i^*\we x_j^*+x_p^*\we x_j^*)=(x_i^*+x_j^*+s_{\gen{i,j}}).t_{pi}\end{equation*}
which in turn gives the equality
\begin{equation*} x_i^*+x_j^*+s_{\gen{i,j}}+x_p^*+x_j^*+s_{\gen{p,j}}=x_i^*+x_p^*+x_j^*+s_{\gen{i,j}}+s_{\gen{p,j}},\end{equation*}
hence $x_j^*=0$, a contradiction. 

On the evidence we have so far we propose 

{\bf Conjecture}. Suppose $G=\displaystyle\prod_{i=1}^mC_{p^{e_i}}^{n_i},\,e_1<\cdots e_m$. Let $N(G,p)$ be the number of almost abelian Hopf algebras of dimension $|G|p$. The function $N(G,p)$ is a polynomial over $\mathbb{Z}$ of degree $\le e_m$ for all $p\ge e_1+\cdots +e_m$.

\end{document}